\documentclass{article}
\usepackage{amsthm,amsmath,epsf,amssymb,xy}
\xyoption{all}

\newcommand{\cA}{{\mathcal A}}

\newcommand{\dixx}{{\Delta_{X \times X}}}

\newcommand{\Da}{{\mathcal D}}
\newcommand{\Der}{{\bf {\rm D} }}

\newcommand{\Du}{{\mathcal D}}
\newcommand{\Ea}{{\mathcal E}}

\newcommand{\tLE}{{\mathcal L \mathcal E}}
\newcommand{\tRE}{{\mathcal R \mathcal E}}
\newcommand{\tEL}{{\mathcal E}_L}
\newcommand{\tER}{{\mathcal E}_R}

\newcommand{\tLEF}{{\mathcal L (\mathcal E, \mathcal F)}}
\newcommand{\tREF}{{\mathcal R (\mathcal E, \mathcal F)}}

\newcommand{\tG}{\mathcal G}
\newcommand{\tF}{\mathcal F}

\newcommand{\REF}{{\mathcal R^{\prime}(\mathcal E, \mathcal F)}}

\newcommand{\tLEE}{{\mathcal L (\mathcal E, \mathcal E)}}
\newcommand{\tREE}{{\mathcal R (\mathcal E, \mathcal E)}}
\newcommand{\tLFE}{{\mathcal L (\mathcal F, \mathcal E)}}
\newcommand{\tRFE}{{\mathcal R (\mathcal F, \mathcal E)}}

\newcommand{\kb}{\bar{k}}

\newcommand{\Fa}{{\mathcal F}}

\newcommand{\Ga}{{\mathcal G}}
\newcommand{\Ha}{{\mathcal H}}
\newcommand{\Hom}{{ {\rm Hom} }}

\newcommand{\Ka}{{\mathcal K}}

\newcommand{\La}{{\mathcal L}}

\newcommand{\Lb}{{\mathbf L}}
\newcommand{\Ltensor}{\mathbin{\overset{\mathbf L}\otimes}}

\newcommand{\Oa}{{\mathcal O}}

\newcommand{\PP}{{\mathbb P}}

\newcommand{\Rb}{{\mathbf R}}

\newcommand{\Rhh}{{\mathbf R {\mathcal Hom}}}

\newcommand{\Za}{{\mathcal Z}}
\newcommand{\wt}[1]{{\widetilde{#1}}}
\newcommand{\sito}{{\overset{\cong}{\longrightarrow}}}

\newtheorem{thm}{Theorem}[section]

\newtheorem{definition}[thm]{Definition}
\newtheorem{proposition}[thm]{Proposition}
\newtheorem{lemma}[thm]{Lemma}

\theoremstyle{definition}
\newtheorem{remark}[thm]{Remark}

\newtheorem{example}[thm]{Example}

\newif\iffigs\figstrue

\numberwithin{equation}{section}

\begin{document}

\begin{center}

{\Large Derived Category Automorphisms from \\Mirror Symmetry}

\end{center}
\bigskip
\centerline{R. Paul Horja}
\bigskip
\centerline{\it Department of Mathematics}
\centerline{\it University of Michigan}
\centerline{\it Ann Arbor, MI 48109, USA}
\centerline{\tt horja@umich.edu}

\begin{abstract} Inspired by the homological mirror symmetry conjecture of
Kontsevich \cite{Kont1},
we construct new classes of automorphisms of
the bounded derived category of coherent sheaves on a smooth quasi--projective
variety. {\it MSC (2000): 18E30; 14J32.}
\end{abstract}

\section{Introduction}
\label{cha:introd}

Let $X$ be a quasi--projective variety over an algebraically closed 
field $\kb,$ $E$ a {\it proper} subvariety of $X$ of codimension $d$
(with $i$ denoting the closed embedding $E \hookrightarrow X$), 
and $q: E \to Z$ a {\it flat} morphism. We assume that $X,E$ and 
$Z$ are 
{\it smooth} varieties
of dimensions $n, n-d$ and $n-d-k$ respectively, as indicated in the 
diagram below.   

\begin{equation}\label{diag:gen}
\begin{split}
\xymatrix{
E^{n-d} \; \ar@{^{(}->}[r]_i  \ar[d]_q & X^n \\
Z^{n-d-k}&
}
\end{split}
\end{equation}
We assume that there exists an invertible sheaf $\theta$ on $Z,$ 
such that $q^* ( \theta \Ltensor \omega_Z^{-1}) \cong \Lb i^* \omega_X^{-1},$
where $\omega_Z$ and $\omega_X$ are the canonical sheaves of $Z$ and $X,$  
respectively. Note that if $\omega_X$ is trivial, then $\theta=\omega_Z$ has the 
required property. 

The goal of this note is to construct new classes of automorphisms of
the bounded derived category of coherent sheaves on $X,$ 
associated with the geometrical context sketched above. An automorphism of a triangulated 
category is an exact functor which induces a self--equivalence. 
The automorphisms defined here are determined by the
$EZ$--spherical objects $\Ea$ (definition \ref{def:EZ})
in the bounded derived category of
coherent sheaves of $E$
that generalize the spherical objects introduced by Kontsevich \cite{Kont2},
and Seidel and Thomas \cite{ST}.
The main example and inspiration for this work is the case
when $X$ is a projective Calabi--Yau variety.

Following the work of Mukai \cite{Mukai1} on actions of functors
on derived categories of coherent
sheaves, Bondal and Orlov \cite{BO}, Orlov \cite{O}
and Bridgeland \cite{Bridge} established criteria that characterize the
equivalences of derived categories of coherent sheaves (see remark 
\ref{rembobr}). Their work
showed that in many instances the group of derived
self--equivalences (automorphisms) captures essential properties of the
algebraic variety itself. Quite tellingly, the case of Calabi--Yau varieties
turned out to be one of the most interesting and difficult to study. Similar
techniques have emerged as useful tools in various other contexts,
for example, in the study of
the McKay correspondence \cite{BKR}, birational geometry \cite{Bridge1},
heterotic string theory compactifications \cite{AspDon} and
D--branes in string theory \cite{Douglas1}, \cite{PA}.

The foundation for our line of investigation is provided by the homological
mirror symmetry conjecture of Kontsevich \cite{Kont1} (see \cite{KontSoib} for
recent progress in this direction). The conjecture states that mirror
symmetry should be viewed an as equivalence between ``Fukaya
$A_\infty$ category'' of a
Calabi--Yau variety $Y$ and the bounded derived category of coherent
sheaves $\Der (X)$ of the mirror Calabi--Yau variety $X.$
The point of view of the present work
stems from some ramifications of the general conjecture as presented by
Kontsevich in \cite{Kont2}. Important results in this direction were obtained
in the paper of Seidel and Thomas \cite{ST}, in which the derived category
braided automorphisms represented by Fourier--Mukai functors associated with
spherical objects were analyzed and compared with the mirror generalized
Dehn twists of Seidel \cite{seidel}. An investigation (adapted to
the toric geometry context) of the correspondence between the expected mirror
automorphisms was pursued by the present author in \cite{rph}.

In fact, the classes of automorphisms introduced in the present work
were obtained with the significant guidance
provided by the interpretation of the
mirror symmetric calculations presented in \cite{rph}.
As a general principle, the automorphisms of the bounded derived category of
coherent sheaves on a Calabi--Yau variety $X$ should be be mirrored by some
automorphisms of the Fukaya category of the mirror Calabi--Yau variety $Y.$
Some of the latter automorphisms are determined by loops in the moduli
space of complex structures on $Y$ and they will depend on the type
of components of the discriminant locus in the moduli space of complex structures
on $Y$ that are surrounded by the loop.
Our proposal
states that for each component of the discriminant
locus in the moduli space of complex structures on $Y$ there is
{\em a whole class} of automorphisms of the bounded derived category of coherent
sheaves on $X$ induced by the $EZ$--spherical objects associated with a
diagram of the type (\ref{diag:gen}). For example, in the toric case,
the so--called ``A--discriminantal hypersurface'' of \cite{GKZ2}
(also called the ``principal component'' in \cite{CK}, \cite{rph} or the
``conifold locus'' in physics)\footnote
{Following a suggestion of Paul Aspinwall \cite{PA}, we choose
to call it the {\it primary} component of the discriminant locus.}
in the moduli space of
complex structures on $Y$ corresponds to the class of spherical objects
on $X$ (in the sense of Seidel and Thomas) which is obtained
in our context for $Z={\rm Spec} (\kb)$ (example \ref{exam:st}). At least 
in the Calabi--Yau case, 
the more general $EZ$--spherical objects introduced in the present work
and their associated automorphisms are geometric manifestations
of the so--called ``phase transitions'' in string theory  \cite{Witten1},
\cite{AGM1}. Elementary contractions in the sense of Mori theory
are examples of these (see example \ref{exam:mori})
and, through mirror symmetry, they
will correspond to the various other components of the discriminant locus in
the moduli space of complex structures on $Y.$
A precise statement
about how this correspondence should work in the toric case was given
in \cite{rph}. An analysis of the correspondence in the closely related language
of D--branes in string theory has been recently presented by P.~Aspinwall \cite{PA}.

There are further questions that can be posed in this context.
The structure of the group of automorphisms is quite intricate and in general
difficult to handle.
Nevertheless, it seems that the structure
of the discriminant locus in moduli space of complex structures on $Y$
encodes a lot of information about the group of automorphisms
of the bounded derived category of coherent
sheaves on $X.$ As explained to me by Professors V.~Lunts and Y.~Manin
(see the end of \cite{maninecm}),
the motivic version of the group should be intimately related to the
Lie algebra actions studied by Looijenga and Lunts \cite{LoLu}.
The relationship has been established in the case of abelian varieties
in \cite{glo}. A related proposal has been recently analyzed
by B.~Szendr\H{o}i \cite{szen}.

Another facet of the story concerns the mirror symmetric Fukaya
category automorphisms. Are there any ``$EZ$--generalized Dehn twists''
associated to (special) Lagrangian vanishing cycles in $Y$
that are not necessarily topological spheres $S^n$?

The automorphisms that are introduced in this work are local in the
sense that they twist only the ``part'' of
the derived category $\Der (X)$ that consists of objects
supported on the subvariety $E.$  The property of $\Ea \in \Der (E)$
to be $EZ$--spherical for a configuration of the type shown in the diagram
(\ref{diag:gen}) does not depend on the ambient space $X$ (with $\omega_X=0$),
but rather on the normal bundle $N_{E/X}$ (remark \ref{rem:normal}).
The local character of the picture
has direct links with the so--called
``local mirror symmetry'' frequently invoked by physicists \cite{KatzMayrVafa},
\cite{CKYZ}. By using the theory of
$t$--structures on triangulated categories \cite{BBD}, this rough idea of ``locality''
can be made more precise. There is now a related and
very interesting proposal in physics made by Douglas \cite{Douglas1} that 
inspired a remarkable construction in mathematics by Bridgeland \cite{Bridgenew} 
about how to express the stability of D--branes in string theory using the triangulated
structure of derived categories.

{\bf Acknowledgments.} It is a pleasure to thank P.~Aspinwall,
D.~Morrison,
P.~Seidel, B.~Szendr\H{o}i and R.~Thomas for very useful comments
which helped improve this work. 
Of course, the 
responsibility for all the statements made in this work rests with me.
Some of the results presented here were obtained
during my stay at the Max--Planck--Institut f\"ur Mathematik in Bonn,
Germany. I am grateful to MPIM and IAS for support and hospitality, and
to Professor Y.~Manin for understanding. The author was supported at the
IAS, Princeton by the NSF grant DMS 97-29992.

{\it Addendum.} This paper is (hopefully) an improved version of an
earlier preprint. I am very grateful to the careful referees who found a gap in the 
earlier version that led to a major revision. 
Since the writing of the preprint, the derived category techniques 
have proved to be very useful in many situations in algebraic geometry and 
string theory. We make no attempt to relate this work to these developments,
although there seem to be many interesting connections. We only mention
the papers \cite{balazs1}, \cite{balazs2}, \cite{alishii}, \cite{ahk},
which are intimately related to the results of the present work and provide
more applications than those presented in the last section of this paper.

\

\subsection{Notation and other general considerations}
\label{cha:not}

A variety is an integral separated scheme of finite type
over an algebraically closed field $\kb.$ 
All the schemes considered in this work are quasi--projective
and Gorenstein (section V.9 in \cite{hart1}), i.e.
the structure sheaf is a dualizing sheaf. For a Gorenstein scheme $W$ over 
$\kb,$ $\Gamma : W \to {\rm Spec} (\kb),$ of 
pure dimension $w,$ the dualizing sheaf $\omega_W$ is invertible and has
the property $\Gamma^!(k)=\omega_W [w].$ Proposition V.9.6 in 
\cite{hart1} shows that, if $E$ and $Z$ are Gorenstein schemes of pure dimensions 
$n-d$ and $n-d-k,$ respectively, and 
$q:E \to Z$ is flat, then  
the scheme $E \times_Z E$ is Gorenstein of pure dimension $n-d+k.$

For a quasi--projective scheme $W,$ we denote by $\Der  (W)$ 
the bounded
derived category of coherent sheaves on $W$
(see \cite{verdier}, \cite{deligne1}, \cite{hart1}, \cite{gelmanin}, 
\cite{BConrad}, for the general
theory of derived and triangulated categories and duality theory).

We use the commonly accepted notation for the derived functors between derived
categories, and
the notation $[n]$ for the shift by $n$ functor in a triangulated
category. The Verdier (derived) dual of an object $\Fa \in \Der (W)$ will be
denoted by $\Du_W \Fa,$ with $ \Du_W \Fa \in \Der (W)$ defined by
\begin{equation*}
\Du_W \Fa := \Rhh_W (\Fa, \Oa_W).
\end{equation*}

For a proper morphism of noetherian schemes of finite Krull dimension
$f : W \to V,$ the construction of a
right adjoint $f^!$ for $\Rb f_*: \Der (W) \to \Der (V)$ was originally established
in the algebraic context by Grothendieck (see \cite{hart1}, \cite{deligne1},
\cite{verdier0}, \cite{BConrad}). We will use  extensively the properties 
associated with the
triple of adjoint functors $(\Lb f^*, \Rb f_*, f^!).$ Given
$\Fa \in \Der (W),$ $\Ga \in \Der (V),$ we will also use the natural
functorial isomorphisms
\begin{equation}
\label{eq:iso1111}
\Rb f_* \Rhh_W ( \Lb f^* \Ga, \Fa) \cong \Rhh_V ( \Ga, \Rb f_* \Fa)
\qquad \text{(\cite{hart1} II.5.10)}
\end{equation}
and
\begin{equation}
\label{eq:dual}
\Rb f_* \Rhh_W ( \Fa, f^! \Ga) \cong \Rhh_V ( \Rb f_* \Fa, \Ga),
\end{equation}
with the latter one provided by the Grothendieck duality theorem (\cite{hart1} III.11.1).

In the concrete situation presented in the beginning of the introduction
(diagram (\ref{diag:gen})), and under the given hypotheses,
the functors $i^!$ and $q^!$ are given by
\begin{equation*}\label{eq:!morph}
\begin{split}
i^! (-) &\cong \Lb i^* (-) \Ltensor i^! ( \Oa_X) 
\cong \Lb i^* (-) \Ltensor (\omega_E \Ltensor \Lb i^* \omega_X^{-1}) [-d],\\
q^! (-) &\cong q^* (-) \Ltensor q^! ( \Oa_Z) = q^* (-) \Ltensor \omega_{E/Z} [k]
\cong q^* (-) \Ltensor (\omega_E \Ltensor q^* \omega_Z^{-1}) [k].
\end{split}
\end{equation*}

\section{Main Results}

We work under the assumptions presented in the beginning of the introduction.
Our particular geometric context leads us to consider the following commutative diagram
($Y:= E \times_Z E$). One could say that, in some sense, this note is
about understanding the geometry summarized by this diagram.

\begin{equation}\label{diagr:big}
\begin{split}
\xymatrix{
  &   & X \times X \ar@/_1pc/[llddd]_{p_2} \ar@/^1pc/[rrddd]^{p_1}                   & & \\
  &   & E \times E \ar@{^{(}->}[u]^l \ar[ddl]_{r_2} \ar[ddr]^{r_1}     & & \\
  &   & \overset{\;}Y \ar@{^{(}->}[u]^k \ar[dl]^{q_2} \ar[dr]_{q_1}
\ar@/_2pc/[uu]_j |!{[u];[dr]}\hole
\ar[dd]^t  & &\\
  X &\;E \ar@{_{(}->}[l]^{i} \ar[dr]^q &
  & E\; \ar@{^{(}->}[r]_{i} \ar[dl]_q &X\\
  &                         & Z &                         &}
\end{split}
\end{equation}
We assume that 
there exists an invertible sheaf $\theta$ on $Z$ such that 
\begin{equation}\label{def:theta}
q^! \theta [-d-k] \cong i^! \Oa_X \cong \omega_{E/X} [-d],
\end{equation}
which is the same thing as 
\begin{equation}\label{def:theta1}
q^* ( \theta \Ltensor \omega_Z^{-1}) \cong \Lb i^* \omega_X^{-1}.
\end{equation}
In other words, the condition says that the pull--back of
$\omega_X$ to $E$ is also the pull-back of an invertible 
sheaf on $Z.$
If $\Rb q_* (\Oa_E) \cong \Oa_Z,$ and $\theta$
has the above property, 
the projection formula implies that $\theta \cong \omega_Z \Ltensor
\Rb q_* (\Lb i^* \omega_X^{-1}).$


\subsection{$EZ$--spherical objects}
\label{cha:ezdef}

Since $i: E  \hookrightarrow X$ is a regular embedding of codimension $d,$ 
$d \geq 0,$ the normal sheaf of $E$ in $X,$ 
$\nu = ({\mathcal J}_E /  {\mathcal J}_E^2)^\vee,$ is locally free, and
$\Lambda^d \nu [-d] \cong \omega_{E/X} [-d] \cong i^! \Oa_X.$

\begin{definition}\label{def:EZ}
An object $\Ea \in \Der (E)$ 
is said to be {\bf $EZ$--spherical} if, either, 
\begin{itemize}\item[i)]
$d=0$ 
and there exists a distinguished triangle 
\begin{equation}\label{eq:ttrr}
\Oa_Z 
\longrightarrow
\Rb q_* \Rhh_E \big( \Ea, \Ea \big) 
\longrightarrow
\theta [-k]
\to \Oa_Z [1],
\end{equation}
or, 
\item[ii)]
$d >0,$ and 
\begin{equation}\label{eq:tttrrr}
\Rb q_* \Rhh_E \big( \Ea, \Ea) \cong \Oa_Z, \ 
\Rb q_* \Rhh_E \big( \Ea, \Ea \Ltensor \Lambda^c \nu \big) = 0, 
\ \text{for} \ 0<c<d.
\end{equation}
\end{itemize}
\end{definition}

\begin{remark}\label{rem:mmm}
Note that, in the case $d >0,$ the first condition in $ii)$ implies that 
\begin{equation}\label{eq:tttrrr1}
\Rb q_* \Rhh_E \big( \Ea, \Ea \Ltensor \Lambda^d \nu \big) \cong \theta [-k].
\end{equation}
Indeed, if
$$
\Rb q_* \Rhh_E \big( \Ea, \Ea) \cong \Oa_Z,
$$
then, Grothendieck duality implies that
\begin{equation}\begin{split}\label{eq:grvv}
&\Rb q_* \Rhh_E \big( \Ea, \Ea \Ltensor \Lambda^d \nu \big) \cong
\Rb q_* \Rhh_E \big( \Ea, \Ea \Ltensor \omega_{E/X} \big) \cong\\
&\cong \Rb q_* \Rhh_E \big( \Rhh_E(\Ea, \Ea), q^! \theta [-k] \big) \cong
\Rhh_Z \big( \Rb q_* \Rhh_E(\Ea, \Ea), \theta [-k] \big) \cong 
 \theta [-k]. 
\end{split}
\end{equation}
\end{remark}

\begin{remark} Note also that in the case $d=0,$ an object $\Ea \in \Der(X)$ such that
$$
\Rb q_* \Rhh_E \big( \Ea, \Ea) \cong \Oa_Z \oplus \theta [-k]
$$
is $EZ$--spherical.
If $Z={\rm Spec} (k),$ 
an $EZ$--spherical object is spherical in the sense of
\cite{ST}. Note also the similarity with the notion of a simple
morphism discussed there.
\end{remark}

\subsection
{Kernels and  Derived Correspondences}

Given two schemes $X_1, X_2,$ an object
$\tG \in \Der (X_1 \times X_2),$ with finite Tor-dimension and 
proper support over both factors, determines an exact functor
$\Phi_{\tG} : \Der (X_1) \to \Der (X_2)$ by the formula (\cite{Mukai1} prop. 1.4)
\begin{equation*}
\Phi_{\tG}(?):=\Rb p_{2_*}
(\tG \Ltensor p_1^*(?)),
\end{equation*}
The object $\tG \in \Der (X_1 \times X_2)$ is called
a {\it kernel}. By analogy with the classical theory of correspondences
(\cite{fulton} chap. 16), a kernel in $\Der (X_1 \times X_2)$ is also said to be
a {\it correspondence} from $X_1$ to $X_2.$

In the situation described by diagram (\ref{diagr:big}),  
$E$ is proper,
so any object $\Ga^{\prime} \in \Der (E\times E)$ 
induces an exact functor
$\Phi_{\tG} : \Der(X) \to \Der(X)$ by 
\begin{equation*}
\Phi_{\tG} (?) := \Rb p_{2_*} (\tG \Ltensor p_1^*(?)), 
\end{equation*}
where $\tG := l_* (\Ga^{\prime}) \in \Der(X \times X).$ 

By the projection formula (II.5.6 in \cite{hart1}), we have that
\begin{equation*}
i_* \Phi_{\Ga^{\prime}} \Lb i^* (?) \cong i_* \Rb r_{2_*} 
(\Ga^{\prime} \Ltensor r_1^* \Lb i^* (?))
\cong \Rb p_{2_*} (\tG \Ltensor p_1^*(?)) \cong \Phi_{\tG} (?).
\end{equation*}

We recall some well--known facts about kernels (correspondences) between 
bounded derived categories of coherent sheaves and their induced 
exact/Fourier--Mukai functors in the case of smooth quasi--projective
varieties.

All the kernels that appear in this work have proper support over the two factors.
The composition of the kernels (correspondences) $\tG_{12} \in \Der (X_1 \times X_2)$
and $\tG_{23} \in \Der (X_2 \times X_3)$ is defined by
\begin{equation} \label{def:comp}
\tG_{23} \star \tG_{12} :=\Rb p_{13_*} \big( p_{12}^* (\tG_{12}) \Ltensor
p_{23}^* (\tG_{23}) \big),
\end{equation}
with $\tG_{23} \star \tG_{12} \in \Der  (X_1 \times X_3).$ 
In the case $X_1=X_2=X,$ the neutral element for the composition of
kernels in $\Der  (X \times X)$ is $\Delta_* \Oa_X,$ where,
as before, $\Delta : X \hookrightarrow X \times X$ is the diagonal morphism.
It is well known (\cite{Mukai1} prop. 1.3.) that there is a natural isomorphism
of functors
\begin{equation*}
\Phi_{\tG_{23} \star \tG_{12}}(?) \cong \Phi_{\tG_{23}} \circ \Phi_{\tG_{12}}(?).
\end{equation*}

From the formula (\ref{def:comp}), it follows that the functors
$\Theta_{\tG_{23}}: \Der (X_1 \times X_2)\to \Der (X_1 \times X_3)$ and 
$\Theta_{\tG_{12}}: \Der (X_2 \times X_3) \to \Der (X_1 \times X_3)$) 
defined by
\begin{equation*}
\begin{split}
\Theta_{\tG_{23}} (?):= \tG_{23} \star (?), 
\Theta_{\tG_{12}} (?):= (?) \star \tG_{12},
\end{split}
\end{equation*}
are exact functor between triangulated categories 
and an argument similar
to the proof of proposition 16.1.1. in \cite{fulton} gives the associativity
of the composition of kernels
\begin{equation}\label{eq:assoc}
\begin{split}
&\tG_{34} \star (\tG_{23} \star \tG_{12}) \cong (\tG_{34} \star \tG_{23}) \star \tG_{12} \cong \\
&\cong \Rb (p^{1234}_{14})_* \Big( 
(p^{1234}_{12})^* (\tG_{12}) 
\Ltensor
(p^{1234}_{23})^* (\tG_{23}) 
\Ltensor
(p^{1234}_{34})^* (\tG_{34})
\Big)
\end{split}
\end{equation}
with $\tG_{12} \in \Der (X_1 \times X_2),$
$\tG_{23} \in \Der (X_2 \times X_3),$ $\tG_{34} \in \Der (X_3 \times X_4).$

\

We now return to the geometric context of the diagram (\ref{diagr:big}).
For any objects $\Ea , \Fa \in \Der (E),$ 
we introduce the following objects
in $\Der  (X \times X)$
\begin{equation}\label{eq:deftrtl}
\begin{split}
\tREF&:= j_*\big(q_1^* (\Du_E \Ea \Ltensor i^! \Oa_X)
\; \Ltensor \; q_2^* (\Fa) \big), \quad \tRE:=\tREE,
\\
\tLEF&:= j_*\big(q_1^* (\Du_E \Ea \Ltensor q^! \Oa_Z)
\; \Ltensor \; q_2^* (\Fa) \big), \quad \tLE:=\tLEE,
\end{split}
\end{equation}
where
\begin{equation*}
\begin{split}
i^! \Oa_X &= 
\omega_E \Ltensor \Lb i^* \omega_X [-d], \\
q^! \Oa_Z &= \omega_{E/Z} [k] = \omega_E \Ltensor q^* \omega_Z^{-1} [k].
\end{split}
\end{equation*}

\begin{proposition} \label{prop:adjoint}
For any $\Ea, \Fa \in \Der (E),$ 
there exist functorial isomorphisms
\begin{equation*}
\begin{split}
i)  \qquad \Hom_{X \times X} (\; ? \star \tREF, \; ??) &\cong 
\Hom_{X \times X} (\;?, \;?? \star \tLFE), \\
\Hom_{X \times X} (\tLEF \star \;?, \; ??) &\cong 
\Hom_{X \times X} (\;?, \tRFE \star \; ??). 
\end{split}
\end{equation*}

\begin{equation*}
\begin{split}
ii)  \qquad \Hom_{X \times X} (\; ? \star \tLEF, \; ??) &\cong 
\Hom_{X \times X} (\;?, \;?? \star \tRFE), \\
\Hom_{X \times X} (\tREF \star \;?, \; ??) &\cong 
\Hom_{X \times X} (\;?, \tLFE \star \; ??). 
\end{split}
\end{equation*}
\end{proposition}

\begin{proof}
This result is a version of a well known fact in the projective case 
(see for example Lemma 1.2 in \cite{BO}). 
In order to be able to use the adjunction theorems in our context, the lack
of properness requires more care. We first prove that
\begin{equation*}
\Hom_{X_1 \times X_3} (\; ? \star \tREF, \; ?? ) \cong 
\Hom_{X_2 \times X_3} (\;?, \;?? \star \tLFE),
\end{equation*}
where the labels for the different copies of $X$ 
will help with the bookkeeping.

We first note that
\begin{equation*} 
? \star \tREF \cong \Rb p_{13_*} \big( 
p_{12}^* (j_{12_*}(\REF) 
\Ltensor
p_{23}^* (?) 
\big),
\end{equation*}
with the appropriate object $\REF \in \Der (Y_{12}),$ 
$\REF= q_1^* (\Du_{E_1} \Ea \Ltensor i^! \Oa_{X_1})
\; \Ltensor \; q_2^* (\Fa),$ $Y_{12}=E_1 \times_Z E_2.$

The fiber square 
\begin{equation}
\begin{split}
\xymatrix{
Y_{12} \times X_3 \ar[r]^{m_{12}} \ar[d]_{(j_{12},Id_{X_3})} & Y_{12} \ar[d]^{j_{12}} 
\\
X_1 \times X_2 \times X_3  \ar[r]^{p_{12}}               & X_1 \times X_2
}
\end{split}
\end{equation}
shows that
\begin{equation*}
p_{12}^* \; j_{12_*} \cong (j_{12},Id_{X_3} )_* \; m_{12}^*,
\end{equation*}
where $m_{12}$ is the projection map $Y_{12} \times X_3 \to Y_{12}.$
By the projection formula, we then obtain that
\begin{equation}\label{eq:kuk}
\begin{split}
\Rb p_{13_*} \big( 
p_{12}^* (j_{12_*}(\REF) 
\Ltensor
p_{23}^* (?) 
\big) &\cong 
\Rb p_{13_*} \; (j_{12},Id_{X_3})_* \big( m_{12}^*(\REF) \Ltensor
(j_{12},Id_{X_3})^* \; p_{23}^* (?) \big) \cong \\
&\cong \Rb k_{13_*} \big( m_{12}^* (\REF) \Ltensor
k_{23}^* (?) \big), 
\end{split}
\end{equation}
where the {\it proper} maps $k_{13}: Y_{12} \times X_3
\to X_1 \times X_3$ and $k_{23}: Y_{12} \times X_3
\to X_2 \times X_3$ are defined as the compositions  
\begin{equation*}
k_{13}:=p_{13} (j_{12}, Id_{X_3}), k_{23}:=p_{23} (j_{12}, Id_{X_3}).
\end{equation*}
By adjunction, we obtain that
\begin{equation}\label{eq:seq100}
\begin{split}
\Hom_{X_1 \times X_3} \big(\; ? \star \tREF, \; ?? ) &\cong 
\Hom_{X_1 \times X_3} \big(
\Rb k_{13_*} \big( m_{12}^* (\REF) \Ltensor
k_{23}^* (?) \big), \; ?? \big) \cong \\
&\cong \Hom_{Y_{12} \times X_3} 
\big( m_{12}^* (\REF) \Ltensor
k_{23}^* (?), k_{13}^! (??) \big) \cong \\
&\cong \Hom_{Y_{12} \times X_3} 
\big(k_{23}^* (?), 
m_{12}^* (\Du_{Y_{12}}\REF) \Ltensor k_{13}^! (??) \big) \cong \\
&\cong \Hom_{X_2 \times X_3} \big(\; ?, \Rb k_{23_*}
\big(m_{12}^* (\Du_{Y_{12}}\REF) \Ltensor k_{13}^! (??) \big)\big).
\end{split}
\end{equation}
We have also used the fact that, since 
the projections $Y_{12} \to E_1,$ $Y_{12} \to E_2$ 
are flat, and $\Ea \in \Der (E_1),$ $\Fa \in \Der (E_2),$
the object $\REF \in \Der (Y_{12})$ is isomorphic to a 
bounded complex of locally free sheaves of finite rank 
(even though $Y_{12}$ might be singular).
The fiber square 
\begin{equation}
\begin{split}
\xymatrix{
Y_{12} \times X_3 \ar[r]^{m_{12}} \ar[d]_{k_{13}} & Y_{12} \ar[d]^{u_1} 
\\
X_1 \times X_3  \ar[r]^{p_1}               & X_1
}
\end{split}
\end{equation}
with $p_1$ flat and $u_1$ proper, and the base change theorem 
(theorem 2 in \cite{verdier0},
or theorem 5(i) in \cite{kleiman}) shows that
\begin{equation*}
m_{12}^* u_1^! \Oa_{X_1} \cong k_{13}^! p_1^* \Oa_{X_1} \cong
k_{13}^!\Oa_{X_1 \times X_3}.
\end{equation*}
By combining this fact with (\ref{eq:seq100}) 
and the adapted version of (\ref{eq:kuk}), 
we see that
\begin{equation}\label{eq:sqe2}
\Hom_{X_1 \times X_3} \big(\; ? \star \tREF, \; ??) \cong 
\Hom_{X_2 \times X_3} \big(?, \Rb p_{23_*}     
(p_{12}^* \big( j_{12_*}(\Du_{Y_{12}}\REF \Ltensor u_1^! \Oa_X) \big)
\Ltensor p_{13}^* (??))\big).
\end{equation}
Recall that $\REF= q_1^* (\Du_E \Ea \Ltensor i^! \Oa_{X})
\; \Ltensor \; q_2^* (\Fa).$ Hence, we can write that
\begin{equation}\label{eq:sqe1}
\begin{split}
\Du_{Y_{12}} \REF \Ltensor u_1^! \Oa_X &\cong
\Rhh_{Y_{12}} (q_1^* (\Du_E \Ea \Ltensor i^! \Oa_X)
\; \Ltensor \; q_2^* (\Fa), \Oa_{Y_{12}}) \Ltensor u_1^! \Oa_{X_1} 
\cong \\
&\cong q_1^* (\Ea) \Ltensor q_2^* (\Du_E \Fa \Ltensor q^! \Oa_Z), 
\end{split}
\end{equation}
since 
\begin{equation}\label{eq:late}
\Du_{Y_{12}} \big(q_1^* (i^! \Oa_X)\big) \Ltensor 
u_1^! \Oa_{X_1} \cong  q_1^* (\omega_E^{-1} \Ltensor
\Lb i^* \omega_X [-d]) \Ltensor (\omega_{Y_{12}} \Ltensor q_1^*\Lb i^* \omega_X^{-1} [-d+k])
\cong 
q_2^* (q^! \Oa_Z),
\end{equation}
where we use proposition II 5.8 in \cite{hart1} ($q_1$ is flat)
and 
the flat base change theorem quoted above.
Combine this fact with (\ref{eq:sqe2}), (\ref{eq:sqe1}), to see that 
we have proved indeed that
\begin{equation*}
\begin{split}
\Hom_{X_1 \times X_3} \big(\; ? \star \tREF, \; ??) &\cong 
\Hom_{X_2 \times X_3} \big(\; ?, \Rb p_{23_*}     
(p_{21}^* \tLFE  \Ltensor p_{13}^* (??))\big) \cong \\
&\cong \Hom_{X_2 \times X_3} \big(\; ?, \; ?? \star \tLFE \big).
\end{split}
\end{equation*}
The proof of the second isomorphism of part $i)$ is completely analogous.

The proof of the isomorphisms of part $ii)$ is also very similar, but  
the existence of the invertible sheaf 
$\theta$ on $Z$ enters the argument in a 
crucial way. For the first isomorphism of part $ii),$
the only change occurs in
formula (\ref{eq:late}) as follows:
\begin{equation*}\begin{split}
\Du_{Y_{12}} \big(q_1^* (q^! \Oa_Z)\big) \Ltensor 
u_1^! \Oa_{X_1} &\cong  
q_2^* \omega_E [k] \Ltensor q_1^* \Lb i^* \omega_X^{-1} [-d+k] \cong 
q_2^* \omega_E \Ltensor t^* ( \theta \Ltensor \omega_Z^{-1}) [-d] \cong \\
&\cong q_2^* \big( \omega_E \Ltensor q^* ( \theta \Ltensor \omega_Z^{-1}) [-d]\big) \cong
q_2^* (i^! \Oa_X).
\end{split}
\end{equation*}
The proof of the second isomorphism of part $ii)$ is again analogous.
\end{proof}

The following easy result will be essential for our arguments.

\begin{lemma}\label{lemma:usus} 

Assume that $\Ka \in \Da$ is an object in a triangulated category 
$\Da$ such that there exists a collection of distinguished triangles
(a ``Postnikov system'') 
\begin{equation}\label{diagr:postnikov}
\begin{split}
\xymatrix@C=9pt{
0=\Ka^{a-1} \ar[rr] & & \Ka^a  \ar[rr] \ar[ld] && \Ka^{a+1} \ar[r] \ar[ld] 
& \dots  \ar[r] & \Ka^{b-1} \ar[rr] & & \Ka^b=\Ka \ar[ld] \\
& \Ha^a (\Ka) [-a] \ar@{.>}[lu] && \Ha^{a+1} (\Ka) [-(a+1)] \ar@{.>}[lu] & & & & 
\Ha^b (\Ka) [-b] \ar@{.>}[lu]
}
\end{split}
\end{equation}
with $a \leq b, a,b \in {\mathbb Z}.$ 
\begin{itemize}
\item[i)] If $T : \Da \to \Da^\prime$ is an exact functor (i.e. triangle 
preserving and commuting with the translations), with $\Da^\prime$
another 
triangulated category, such that 
\begin{equation*}\label{eq:conddd}
T(\Ha^c (\Ka))=0, \ \text{for} \ a< c < b,
\end{equation*}
then there exists a distinguished triangle
\begin{equation*}
\begin{split}
\xymatrix{
T(\Ha^a (\Ka) [-a]) \ar[rr] & &  T(\Ka) 
\ar[ld]\\
& T(\Ha^b (\Ka) [-b]) \ar@{.>}[lu] &
}
\end{split}
\end{equation*}
\item[ii)] If $H : \Da \to \cA$ is a cohomological functor 
(i.e. mapping triangles into long exact sequences), with $\cA$ 
an Abelian category, such that
\begin{equation*}\label{eq:condhh}
H(\Ha^c (\Ka)[-c-p])=0, \ \text{for} \ a< c \leq b, p=0,1,
\end{equation*}
then 
\begin{equation*}
H(\Ka) \cong H(\Ha^a (\Ka) [-a]).
\end{equation*}
\end{itemize}
\end{lemma}

The prototype of a Postnikov system associated to an object $\Ka$ of a triangulated 
category $\Da$ endowed with a t-structure (definition 1.3.1 in \cite{BBD}) 
appears as a succession of distinguished triangles of the form
\begin{equation*}
\begin{split}
\xymatrix{
\tau_{<a} \Ka \ar[rr] & & \tau_{\leq a} \Ka \ar[ld]\\
& \Ha^a (\Ka) [-a] \ar@{.>}[lu] &
}
\end{split}
\end{equation*}
where $\tau_{<a}, \tau_{\leq a}$ are the truncation functors, and $\Ha^a$ the 
cohomology functors with values in the heart of the t-structure. 
We will use this construction in what follows for the case of the 
bounded derived category of coherent sheaves endowed with the usual 
t-structure. Some other Postnikov systems to be used in this note are 
provided by the following lemma.

\begin{lemma}\label{lemma:sga} Let $i : T \hookrightarrow S$ be 
a regular embedding of codimension $d,$ with $T$ proper, 
and let $\nu$ be the normal 
sheaf of $T$ in $S.$ 
\begin{itemize}
\item[i)] The exact functor  $\Lb i^* i_* : \Der(T) \to \Der(T)$ is induced by 
the kernel 
$$
\Gamma_{23_*} (\Oa_T) \star \Gamma_{12_*} (\Oa_T) \in \Der (T \times T)
$$
where 
$\Gamma_{12} : T \hookrightarrow T \times S,$ and
$\Gamma_{23} : T \hookrightarrow S \times T,$ the graph (diagonal) 
embeddings 
determined by $i: T \hookrightarrow S,$ and there exists a 
Postnikov system of the type (\ref{diagr:postnikov}) with
$a=-d, b=0,$ and 
$$
\Ka= \Gamma_{23_*} (\Oa_T) \star \Gamma_{12_*} (\Oa_T), \, 
\Ha^{-c} (\Ka) \cong \Delta_* \Lambda^{c} \nu^\vee, 
\,  \text{for $0 \leq c \leq d$},
$$
where $\Delta : T \hookrightarrow T \times T$ is the 
diagonal embedding.
\item[ii)]
For any object $\Ga \in \Der (T)$ there 
exists a Postnikov system of the type (\ref{diagr:postnikov}) with
$a=-d, b=0,$ and 
$$
\Ka= \Lb i^* i_* \Ga, \, \Ha^{-c} (\Ka) \cong \Ga \Ltensor 
\Lambda^c \nu^\vee, \,  \text{for $0 \leq c \leq d$}.
$$
\end{itemize}
\end{lemma}

\begin{proof} 
The exact functor $\Lb i^* i_*$ is induced by 
the kernel 
$$
\Gamma_{23_*} (\Oa_T) \star \Gamma_{12_*} (\Oa_T) \cong
\Rb p_{13_*} \big( 
p_{12}^* \Gamma_{12_*} (\Oa_T)
\Ltensor 
p_{23}^* \Gamma_{23_*} (\Oa_T) \big).
$$
Consider the fiber square (with $p_{12}$ flat)
\begin{equation}
\begin{split}
\xymatrix{
T \times T \ar[r]^{q_1} \ar[d]_{(\Gamma_{12}, Id_T)} & T \ar[d]^{\Gamma_{12}} \\
T \times S \times T \ar[r]^-{p_{12}}               & T \times S
}
\end{split}
\end{equation}
Therefore $p_{12}^* \Gamma_{12_*} (\Oa_T) \cong 
(\Gamma_{12_*}, Id_{T_*}) (\Oa_{T \times T}).$
By the projection formula, we can then write
\begin{equation*}
\Gamma_{23_*} (\Oa_T) \star \Gamma_{12_*} (\Oa_T) \cong
\Rb p_{13_*} (\Gamma_{12_*}, Id_{T_*}) \Big( (\Lb \Gamma_{12}^*, Id_T^*) \;
p_{23}^* \Gamma_{23_*} (\Oa_T) \Big).
\end{equation*}
But $p_{13} \circ (\Gamma_{12}, Id_T) = Id_{T \times T},$ and
$p_{23} \circ (\Gamma_{12}, Id_T) = (i, Id_T).$ Hence
\begin{equation*}
\Gamma_{23_*} (\Oa_T) \star \Gamma_{12_*} (\Oa_T) \cong
(\Lb i^*, Id_T^*) \; \Gamma_{23_*} (\Oa_T).
\end{equation*}
Set $\iota :=(i, Id_T): T \times T \to S \times T.$
Then $\Gamma_{23} = \iota \circ \Delta,$ with $\Delta : T \hookrightarrow
T \times T$ the diagonal embedding, and
it follows that
\begin{equation*}
\Gamma_{23_*} (\Oa_T) \star \Gamma_{12_*} (\Oa_T) 
\cong
\Lb \iota^* \iota_* (\Delta_* \Oa_T). 
\end{equation*}
Now note that there exists a Postnikov system in $\Der (T \times T)$ 
of the type (\ref{diagr:postnikov})  with $\Ka=
\Lb \iota^* \iota_* (\Delta_* \Oa_T),$ and 
$\Ha^{-c}(\Ka)\cong \Lb^{c} \iota^* \iota_* (\Delta_* \Oa_T).$
On $S \times T$ these sheaves are the Tor-sheaves
${\mathcal Tor}_{c}^{\Oa_{S \times T}}
(\iota_* (\Delta_* \Oa_T), \iota_* \Oa_{T \times T} )$ and can be 
computed with the help 
of a local Koszul resolution
of $\iota_* \Oa_T$ in $S \times T$ (see, for example, VII.2.5 in \cite{sga6}). 
It follows that 
$$
\Ha^{-c}(\Ka) \cong \left\{ \begin{array}{ll} \Delta_* \Lambda^{c} \nu^\vee
& \text{if $0 \leq c \leq d,$} \\
0 & \text{otherwise.} \end{array} \right.
$$
The lemma follows, since this Postnikov system on $T \times T$ is 
mapped by an exact functor 
into the required Postnikov system on $T.$ 
\end{proof}

\begin{remark} By combining the previous two lemmas, it follows that, 
for an $EZ$--spherical object $\Ea \in \Der (E),$ there exists a 
distinguished triangle 
$$
\Oa_Z \longrightarrow 
\Rb q_* \Rhh_E \big( \Ea, i^! i_* \Ea) 
\longrightarrow \theta [-d-k] \longrightarrow \Oa_Z [1].
$$
This remark will not be used in this paper.
\end{remark}

\begin{proposition}\label{prop:uuuse} 
If $\Ea$ in $\Der(E)$ is $EZ$--spherical, then: 

i) For $e < d+k,$
\begin{equation*}
\Hom_{X\times X} \big(\tRE, \Delta_* \Oa_X [e]\big) \cong 
\Hom_{X\times X} \big( \Delta_* \Oa_X[-e], \tLE \big) 
\cong \Hom_{Z \times Z} \big(\Delta_* \Oa_Z, \Delta_* \Oa_Z [e] \big).
\end{equation*}
In particular,
\begin{equation*}
\Hom_{X\times X} \big(\tRE, \Delta_* \Oa_X [e] \big)  \cong 
\Hom_{X\times X} \big( \Delta_* \Oa_X [-e], \tLE \big) 
\cong \left\{ \begin{array}{ll} \kb & \text{if $e=0,$} \\
0 & \text{if $e<0.$} \end{array} \right.
\end{equation*}

ii) For $e < 2(d+k),$
\begin{equation*}
\Hom_{X\times X} \big(\tLE, \Delta_* \Oa_X [e]\big) \cong 
\Hom_{X\times X} \big( \Delta_* \Oa_X [-e], \tRE \big) 
\cong \Hom_{Z \times Z} \big( \Delta_* \Oa_Z, \Delta_*
\theta [e-d-k] \big).
\end{equation*}
In particular, for $e < d+k,$
\begin{equation*}
\Hom_{X\times X} \big(\tLE, \Delta_* \Oa_X [e] \big) \cong 
\Hom_{X\times X} \big( \Delta_* \Oa_X [-e], \tRE \big) = 0.
\end{equation*}

iii) For $e < d+k,$
\begin{equation*}
\Hom_{X \times X} \big( \tRE, \tRE [e]\big) \cong
\Hom_{X \times X} \big( \tLE[-e], \tLE \big) \cong 
\Hom_{Z \times Z} \big( \Delta_* \Oa_Z, 
\Delta_* \Oa_Z [e]\big).
\end{equation*}
In particular,
\begin{equation*}
\Hom_{X\times X} (\tRE, \tRE [e] ) \cong 
\Hom_{X \times X} ( \tLE [-e], \tLE ) \cong 
\left\{ \begin{array}{ll} \kb & \text{if $e=0,$} \\
0 & \text{if $e<0.$} \end{array} \right.
\end{equation*}
\end{proposition}
\begin{proof}
Note that the first isomorphisms 
in each of the three parts of proposition 
follow directly from proposition \ref{prop:adjoint}, while the 
isomorphisms in the second groups of each part follow easily
provided we proved the second isomorphisms in each part. 
The notation follows diagram (\ref{diagr:big}).

{\it i)} By adjunction, we have that
\begin{equation}\label{eq:dix0}
\Hom_{X \times X} \big( \tRE, \Delta_* \Oa_X [e] \big) 
\cong \Hom_{E \times E} \Big( k_* \big(q_1^* (\Du_E \Ea \Ltensor i^! \Oa_X)
\; \Ltensor \; q_2^* (\Ea) \big), l^! \Delta_* \Oa_X [e] \Big).
\end{equation}
Since $q_1= r_1 \circ k,$  $q_2=r_2 \circ k,$ the projection formula gives
\begin{equation}
k_* \big(q_1^* (\Du_E \Ea \Ltensor i^! \Oa_X)
\; \Ltensor \; q_2^* (\Ea) \big) \cong  \big(r_1^* (\Du_E \Ea \Ltensor i^! \Oa_X)
\; \Ltensor \; r_2^* (\Ea) \big) \Ltensor k_* \Oa_Y .
\end{equation}
On the other hand
\begin{equation}
l^! \Delta_* \Oa_X \cong \Lb l^* \Delta_* \Oa_X 
\Ltensor l^! \Oa_{X \times X}.
\end{equation}

Note that for $\Lb l^* \Delta_* \Oa_X,$ there exists a Postnikov system 
in $\Der (E \times E)$ 
of the type (\ref{diagr:postnikov}) with $\Ka= \Lb l^* \Delta_* \Oa_X,$ 
and $\Ha^{-c}(\Ka)\cong \Lb^{c} l^* \Delta_* \Oa_X.$ Indeed, 
on $X \times X,$ these sheaves are in fact the Tor-sheaves
${\mathcal Tor}_{c}^{\Oa_{X \times X}}
(\Delta_* \Oa_X, l_* \Oa_{E \times E} ).$ 
Since $l_* \Oa_{E \times E} \cong p_1^* (i_* \Oa_E) \Ltensor p_2^* (i_* \Oa_E)$
we have that 
$$
\Lb \Delta^* l_* \Oa_{E \times E} \cong i_* \Oa_E \Ltensor i_* \Oa_E \cong 
i_* ( \Lb i^* i_* \Oa_E),
$$
so $\Lb^{c} \Delta^* l_* \Oa_{E \times E} \cong i_* \Lambda^{c} \nu^\vee,
$ for $0 \leq c \leq d.$ On $X \times X,$
these sheaves are again 
the Tor-sheaves ${\mathcal Tor}_{c}^{\Oa_{X \times X}}
(\Delta_* \Oa_X, l_* \Oa_{E \times E} ).$ Therefore, we see that
$$
\Ha^{-c}(\Ka) \cong \left\{ \begin{array}{ll} \Delta_* \Lambda^{c} \nu^\vee
& \text{if $0 \leq c \leq d,$} \\
0 & \text{otherwise,} \end{array} \right.
$$
where $\Delta : E \hookrightarrow E \times E$ is the 
diagonal embedding. 

Our goal is to apply part $ii)$ of 
lemma \ref{lemma:usus} for 
the cohomological functor $H : \Der (E) \to {\bf Ab}$ 
induced by the right hand side of (\ref{eq:dix0}),
$$ 
H(\Ka):=  \Hom_{E \times E} \Big( k_* \big(q_1^* (\Du_E \Ea \Ltensor i^! \Oa_X)
\; \Ltensor \; q_2^* (\Ea) \big), \Ka \Ltensor l^! \Oa_{X\times X} [e] \Big).
$$
where $\Ka=\Lb l^* \Delta_* \Oa_X$ with the Postnikov system described above,
and $a=-d, b=0.$ 

The hypotheses of part $ii)$ of lemma \ref{lemma:usus} require us to 
show that the group 
\begin{equation}\label{eq:grp1}
\Hom_{E \times E} \Big( k_* \big(q_1^* (\Du_E \Ea \Ltensor i^! \Oa_X)
\; \Ltensor \; q_2^* (\Ea) \big), \Delta_* \Lambda^{c} \nu^\vee[c-p] 
\Ltensor l^! \Oa_{X \times X}
[e] \Big)
\end{equation}
is zero for $0 \leq c < d, p=0,1,$ and $e < d+k.$

To prove this claim, note 
first that the projection formula implies that
$$
k_* \big(q_1^* (\Du_E \Ea \Ltensor i^! \Oa_X)
\; \Ltensor \; q_2^* (\Ea) \big) \cong 
\big(r_1^* (\Du_E \Ea \Ltensor i^! \Oa_X)
\; \Ltensor \; r_2^* (\Ea) \big) \Ltensor k_* \Oa_Y,
$$
and the already quoted
flat base change theorem (theorem 2 in \cite{verdier0},
or theorem 5(i) in \cite{kleiman}) shows that
\begin{equation}\label{eq:chl}
l^! \Oa_{X \times X} \cong r_1^* i^! \Oa_X \Ltensor r_2^* i^! \Oa_X.
\end{equation}

We can then apply adjunction for the pair of functors 
$\big( \Lb \Delta^*, \Delta_* \big)$ associated to the diagonal 
of $E \times E,$ to obtain that the group (\ref{eq:grp1}) is isomorphic 
to 
\begin{equation*} 
\begin{split}
&\Hom_E  \big( \Lb \Delta^* (k_* \Oa_Y)
\Ltensor (\Du_E \Ea \; \Ltensor \; \Ea),
\Lambda^{c} \nu^\vee[c-p] \Ltensor i^! \Oa_X [e] \big) \cong \\
&\cong \Hom_{E \times E} \big( k_* \Oa_Y, \Delta_*
\Rhh_E (\Ea, \Ea \Ltensor \Lambda^{d-c} \nu \, [-d+c-p] [e] )
\big).
\end{split}
\end{equation*}

Now consider the following
commutative diagram

\begin{equation}\label{diagr:zz}
\begin{split}
\xymatrix@C=4pc{
  &   Y  \ar[r]^t  \ar[d]^k & Z   \ar[d]^{\Delta}  \\
E \ar[ur]^{\Delta_Y} \ar[r]^{\Delta} &   E \times E  \ar[r]^{q \times q}
& Z \times Z }
\end{split}
\end{equation}
where $t = q_1 \circ q = q_2 \circ q.$ Since $q \times q : E \times E \to Z \times Z$
is flat, the base change theorem for the above fiber square allows us to write
\begin{equation}\label{eq:kk}
k_* \Oa_Y \cong k_* t^* \Oa_Z \cong (q \times q)^* \Delta_* (\Oa_Z).
\end{equation}
By adjunction, we conclude that the group (\ref{eq:grp1}) is isomorphic 
to 
\begin{equation}\label{ref:rere}
\Hom_{Z \times Z} \big( \Delta_* \Oa_Z, \Delta_* \Rb q_*
\Rhh_E (\Ea, \Ea \Ltensor \Lambda^{d-c} \nu \, [-d+c-p] [e] )
\big).
\end{equation}
Since $\Ea$ is $EZ$--spherical, we have that 
$\Rhh_E (\Ea, \Ea \Ltensor \Lambda^{d-c} \nu) \cong 0$ for 
$0 < c <d.$ For $c=0,$
remark \ref{rem:mmm} shows that 
$\Rhh_E (\Ea, \Ea \Ltensor \Lambda^{d} \nu)\cong \theta [-k],$
the group above is zero if 
$$
-d -k- p +e < 0,
$$
that is when $e < d+k.$ This proves the above claim that the group 
(\ref{eq:grp1}) is zero for $0 \leq c < d, p=0,1,$ and $e < d+k.$

Therefore, part $ii)$ of 
lemma \ref{lemma:usus} implies indeed that 
$\Hom_{X \times X} \big( \tRE, \Delta_* \Oa_X [e] \big)$
is isomorphic to the group obtained by setting $c=d, p=0$ in (\ref{ref:rere}),
that is 
$$
\Hom_{Z \times Z} \big( \Delta_* \Oa_Z, \Delta_* \Rb q_*
\Rhh_E (\Ea, \Ea [e])  \big), 
$$
which is isomorphic (for an $EZ$--spherical object $\Ea$) to 
$$
\Hom_{Z \times Z} \big( \Delta_* \Oa_Z, 
\Delta_* \Oa_Z [e]
\big).
$$
Note that in the case $d=0,$ the same conclusion is obtained
without having to use lemma \ref{lemma:usus}, since in this case
the $EZ$-condition implies the existence of 
a distinguished triangle of the form 
\begin{equation*}
\begin{split}
\ldots \to \Hom_{Z \times Z} \big( \Delta_* \Oa_Z, \Delta_* \theta [-k+e-1]
\big) &\to
\Hom_{Z \times Z} ( \Delta_* \Oa_Z, 
\Delta_* \Oa_Z [e]
) \to \\
\to \Hom_{Z \times Z} \big( \Delta_* \Oa_Z, \Delta_* \Rb q_*
\Rhh_E (\Ea, \Ea [e])  \big) 
&\cong  \Hom_{X \times X} 
(\tRE, \Delta_* \Oa_X [e]) 
\to \\
\to \Hom_{Z \times Z} \big( \Delta_* \Oa_Z, \Delta_* \theta [-k+e]
\big) &\to \ldots.
\end{split}
\end{equation*}

$ii)$ The proof of the second isomorphism in this part is almost identical
to the previous argument; formula 
(\ref{eq:chl}) has to be replaced by 
\begin{equation*}
l^! \Oa_{X \times X} \cong r_1^* q^! \Oa_Z \Ltensor (r_1^* q^* \theta[-d-k]
\Ltensor r_2^* i^! \Oa_X),
\end{equation*}
since, by (\ref{def:theta}) $q^! \theta [-d-k] \cong i^! \Oa_X.$ 

$iii)$ We now proceed with the proof of the second isomorphism in this part.
We have that
\begin{equation*}\label{eq:dix000}
\Hom_{X \times X} \big( \tRE, \tRE [e]\big) \cong
\Hom_{E \times E} \big( \Za, l^! l_* \Za [e]\big) \cong 
\Hom_{E \times E} \big( \Za, \Lb l^* l_* \Za [e] \Ltensor l^! \Oa_{X \times X}
\big),
\end{equation*}
with $\Za \in \Der (E \times E)$ defined by
\begin{equation*}
\Za:= \big( r_1^* (\Du_E \Ea \Ltensor i^! \Oa_X)
\; \Ltensor \; r_2^* (\Ea) \big) \; \Ltensor \; k_* \Oa_Y.
\end{equation*}

According to lemma \ref{lemma:sga}, there exists a Postnikov 
system with $\Ka=\Lb l^* l_* \Za,$ $a=-2d, b=d$ and
$$
\Ha^{-c}(\Ka) \cong 
\Za \Ltensor \Lambda^{c} 
\tilde{\nu}^\vee,$$ 
for $0 \leq c \leq 2d.$
where $\tilde{\nu}$ is the normal sheaf of $E \times E$ in $X \times X.$

In order to use again part $ii)$ of lemma \ref{lemma:usus} for 
the cohomological functor $H : \Der (E \times E) \to {\bf Ab},$ 
 induced by the 
right hand side of (\ref{eq:dix000}),
$$ 
H(\Ka):=  \Hom_{E \times E} (\Za, \Ka \Ltensor l^! \Oa_{X\times X} [e] \Big).
$$
with $\Ka=\Lb l^* l_* \Za,$ $a=-2d, b=d,$
we need to show that 
\begin{equation}\label{eq:grp2}
\Hom_{E \times E} \big( \Za, \Za [e] \Ltensor  \Lambda^{c} \tilde{\nu}^\vee [c-p]
\Ltensor l^! \Oa_{X \times X} \big),
\end{equation}
is zero for $0 \leq c < 2d, p=0,1,$ and $e < d+k.$ Moreover, the 
K\"unneth formula implies that it is enough to show that the groups
\begin{equation}\label{eq:grp3}
\Hom_{E \times E} \big( \Za, \Za [e] \Ltensor  (\Lambda^{c_1} 
\tilde{\nu_1}^\vee \otimes \Lambda^{c_2}  \tilde{\nu_2}^\vee) [c_1 +c_2-p] 
\Ltensor l^! \Oa_{X \times X} \big),
\end{equation}
are zero for $0 \leq c_1 + c_2 <  2d, p=0,1,$ $e < d+k.$ Here 
the normal sheaves $\tilde{\nu_1}$ and $\tilde{\nu_2}$ give the decomposition 
of the normal sheaf $\tilde{\nu}$ along the two directions corresponding 
to the embeddings $E_i \hookrightarrow X_i.$

After regrouping the various tensor products we obtain the the group 
(\ref{eq:grp3}) is isomorphic to 
\begin{equation*}
\Hom_{E \times E} \Big(  k_* \Oa_Y, 
\big( r_1^* \Rhh_E (\Ea, \Ea \Ltensor \Lambda^{d-c_1} {\nu}
\Ltensor
r_2^* \Rhh_E (\Ea, \Ea \Ltensor \Lambda^{d-c_2} {\nu} \big) 
\Ltensor k_* \Oa_Y [-2d+c_1+c_2-p] [e]\Big).
\end{equation*}
By (\ref{eq:kk}), we have that $k_* \Oa_Y \cong
(q \times q)^* \Delta_* \Oa_Z,$ so
adjunction and the projection formula imply that the group above 
is isomorphic to 
\begin{equation*}
\Hom_{Z \times Z} \Big( \Delta_* \Oa_Z,  \Delta_*
\Rb t_* \big( q_1^* \Rhh_E (\Ea, \Ea \Ltensor \Lambda^{d-c_1} {\nu})
\Ltensor
q_2^* \Rhh_E (\Ea, \Ea \Ltensor \Lambda^{d-c_2} {\nu} )
[-2d+c_1+c_2-p] [e]\Big).
\end{equation*}

Since $q$ is flat, the K\"unneth formula for the interior fiber square
in the diagram (\ref{diagr:big}) says that
\begin{equation*}
\Rb t_* \big( q_1^* \Ga_1 \Ltensor q_2^* \Ga_2 \big)
\cong
\Rb q_* \Ga_1  \Ltensor \Rb q_* \Ga_2,
\end{equation*}
for any $\Ga_1, \Ga_2 \in \Der(E).$ Indeed,
\begin{equation*}
\begin{split}
\Rb q_* \Ga_1  \Ltensor \Rb q_* \Ga_2 &\cong \Rb q_* \big( \Ga_1  \Ltensor q^*
\Rb q_* \Ga_2 \big) \cong \Rb q_* \big( \Ga_1  \Ltensor
\Rb q_{1_*} q_2^* \Ga_2 \big) \\
&\cong \Rb (q \circ q_1)_*
\big( q_1^* \Ga_1 \Ltensor q_2^* \Ga_2 \big).
\end{split}
\end{equation*}
Hence, we have shown that the group (\ref{eq:grp3}) is isomorphic to 
\begin{equation*}
\Hom_{Z \times Z} \Big( \Delta_* \Oa_Z,  \Delta_* \big(
\Rb q_* \Rhh_E (\Ea, \Ea \Ltensor \Lambda^{d-c_1} {\nu})
\Ltensor
\Rb q_* \Rhh_E (\Ea, \Ea \Ltensor \Lambda^{d-c_2} {\nu} ) \big)   
[-2d+c_1+c_2-p] [e]\Big).
\end{equation*}
Clearly, if $\Ea$ is $EZ$--spherical, and $c_1=c_2=0,$ the group is 
\begin{equation*}
\Hom_{Z \times Z} \Big( \Delta_* \Oa_Z,  \Delta_* (\theta \otimes \theta)
[-2d-2k-p] [e]\Big),
\end{equation*}
which clearly is zero for $e < d+k.$ We also have to consider the 
case $c_1+c_2 =d.$ Then the group is 
\begin{equation*}
\Hom_{Z \times Z} \Big( \Delta_* \Oa_Z,  \Delta_* \theta
[-d-k-p] [e]\Big),
\end{equation*}
which is again zero for $e < d+k.$ In all the other situations with 
$0 \leq c_1 + c_2 <2d,$ the group is zero. It is immediate then that part 
$ii)$ of lemma \ref{lemma:usus} implies the required isomorphism
$$
\Hom_{X \times X} \big( \tRE, \tRE [e]\big) \cong 
\Hom_{Z \times Z} \Big( \Delta_* \Oa_Z,  \Delta_* \Oa_Z [e] \Big),
$$
for $e < d+k.$ 

As in the proof of part $i),$ the case $d=0$ from the above calculation
follows without invoking lemma \ref{lemma:usus}.
\end{proof}

For an $EZ$--spherical object $\Ea \in \Der (E),$ 
set $e=0$ in the isomorphisms of part $i)$ of the
above proposition. Denote  
by $r_\Ea$ and $l_\Ea$ the images of the identity morphism
under some choice of isomorphisms at part $i).$
\begin{equation*}
\tRE \overset{r_{\Ea}}\longrightarrow (\dixx)_* (\Oa_X), \quad
(\dixx)_* (\Oa_X) \overset{l_{\Ea}}\longrightarrow \tLE.
\end{equation*}

\begin{definition}\label{def:lr}
If $\Ea \in \Der (E)$ is $EZ$--spherical, 
the objects $\tER$ and $\tEL$ in $\Der (X \times X)$
are defined
(up to non--canonical isomorphisms) by the distinguished
triangles in $\Der(X \times X)$
\begin{equation} \label{eq:basiclr}
\begin{split}
&\tRE  \overset{r_{\Ea}}\longrightarrow
\Delta_* \Oa_X
\to \tER \to \tRE [1] \\
&\tEL \to \Delta_* \Oa_X
\overset{l_{\Ea}}\longrightarrow
\tLE \to \tEL [1].
\end{split}
\end{equation}
\end{definition}

We will need the following corollary of proposition
\ref{prop:uuuse}.

\begin{proposition}\label{prop:corse} 
If $\Ea \in \Der (E)$ is $EZ$--spherical, then 
\begin{equation*}\label{eq:isom400}
\Hom_{X \times X} (\tEL , \tLE \; [-1]) \cong 
\Hom_{X \times X} (\tRE , \tER \; [-1]) = 0.
\end{equation*}
\end{proposition}

\begin{proof}
Apply the cohomological functor $\Hom_{X \times X} (? , \tLE)$
to the distinguished triangle
$$
\tEL \to \Delta_* \Oa_X
\overset{l_{\Ea}}\longrightarrow
\tLE \to \tEL [1],
$$
and look at the following piece of the resulting long exact sequence:
\begin{equation*}
\begin{split}
\ldots &\to \Hom_{X \times X} ( \Delta_* \Oa_X [1], \tLE) 
\to \\
& \to \Hom_{X \times X} (\tEL [1], \tLE) \to
\Hom_{X \times X} (\tLE, \tLE) \overset{l_{\Ea}^{\%}}\longrightarrow
\Hom_{X \times X} ( \Delta_* \Oa_X, \tLE) \to \ldots.
\end{split}
\end{equation*}
Note that, since $\Ea$ is $EZ$--spherical, by $i)$ and $iii)$ of proposition
\ref{prop:uuuse}, the last two groups on the 
right are isomorphic to 
$\Hom_{Z \times Z} \big( \Delta_* \Oa_Z, \Delta_* \Oa_Z \big)$
(and in fact to the field $\kb$). 
The morphism 
$l_{\Ea}$ was chosen to correspond 
to the identity in the group
$\Hom_{Z \times Z} \big(\Delta_* \Oa_Z, \Delta_* \Oa_Z \big),$
which shows that the induced homomorphism $l_{\Ea}^{\%}$
is in fact an isomorphism. Since part $i)$ of the previous proposition
shows that the group
$
\Hom_{X \times X} ( \Delta_* \Oa_X [1], \tLE) 
$  
is zero, we obtain indeed that $\Hom_{X \times X} (\tEL [1], \tLE )$
is zero. The proof of the other half of the proposition is completely
analogous.
\end{proof}

The main result of this note is the following theorem.

\begin{thm} \label{thm} 
Under the hypotheses of definition \ref{def:EZ},
for any $EZ$--spherical object $\Ea \in \Der (E),$ the exact functors
$\Phi_{\tER}$ and $\Phi_{\tEL}$ are inverse automorphisms of $\Der(X),$ i.e.
\begin{equation*}
\Phi_{\tER} \circ \Phi_{\tEL} \cong Id_{\Der(X)}, \qquad
\Phi_{\tEL} \circ \Phi_{\tER} \cong Id_{\Der(X)}.
\end{equation*}
\end{thm}

The proof of the theorem will occupy the next section.

\begin{remark}\label{rembobr}
At least in the projective case, a 
theorem of Bondal and Orlov (theorem 1.1 in
\cite{BO}) and Bridgeland (theorem 5.1 \cite{Bridge})  
provides a clear criterion for the fully faithfulness of the
exact functors $\Phi_{\tER}$ and $\Phi_{\tEL}.$ Namely, $\Phi_{\tER}$ is fully 
faithful if and only if, for each point $x \in X,$ 
$$
\Hom_{X \times X} (\Phi_{\tER}(\Oa_x),\Phi_{\tER}(\Oa_x) )= \kb,
$$
and for each pair of points $x_1, x_2 \in X,$ and for each integer $i,$
$$
\Hom_{X \times X} (\Phi_{\tER}(\Oa_{x_1}),\Phi_{\tER}(\Oa_{x_2})) = 0, \ \text{unless
$x_1=x_2$ and $0 \leq i \leq n.$}
$$
While this is a very geometrical characterization, the quite convoluted 
way of defining the objects $\tER$ and $\tEL$ makes its use rather hard
in the given set-up.
\end{remark}

\begin{remark}\label{rem:goren} As it can be seen from the proofs 
contained in the next section, the smoothness assumptions on
$X, E$ and $Z$ can be relaxed somewhat. Quasi--projectiveness can be replaced 
by assumptions on the schemes 
that guarantee that Grothendieck duality theory works. Beyond that,  
it is enough
to assume that the schemes 
$X,E, Z$ as well as the morphism $q$
are {\it Gorenstein} (\cite{hart1} V.9), with $E$ proper, and that 
$E \hookrightarrow X$ is a regular embedding. 
 We would then have
to add ``by hand'' further assumptions on the kernels $\Ea,$ $\tRE,$
$\tLE;$ for example, we could assume that they are {\it perfect,} namely that 
they are isomorphic (in the corresponding 
derived categories) to bounded complexes of locally free sheaves of finite rank. 
As it can be seen from the proof of proposition \ref{prop:uuuse},
we also need to assume that $Z$ is connected.
The results of this section would remain true.
The crucial assumption that the 
morphism $q$ is proper and flat ensures that the considered functors
take {\em bounded} derived categories of {\em coherent} sheaves to
{\em bounded} derived categories of {\em coherent} sheaves.
\end{remark}

\begin{remark} In the case of a smooth Calabi--Yau variety $X$ ($\omega_X=0,
\theta=\omega_Z$),
the choice $\Ea= \omega_E$ gives
\begin{equation*}
\tRE = j_* ( \Oa_Y [-d] \Ltensor q_2^* \omega_E ).
\end{equation*}
The corresponding kernel $\tER$ is precisely the kernel introduced
in section 4.1 of \cite{rph} in the case of a Calabi--Yau complete
intersection $X$ in a toric variety (that result provided the 
inspiration for this work). 
As it is shown there, under the additional
hypothesis that $E$ is a complete intersection of toric divisors,
the action in cohomology
of the corresponding Fourier--Mukai functor $\Phi_{\tER}$
matches the mirror symmetric monodromy action obtained by
analytical computations. It will be seen in the section \ref{sec:appl} of this work
that, in that context, $\Ea=\omega_E$ is indeed $EZ$--spherical. The
invertability of the corresponding exact functor in the case of a type III birational
contraction for a Calabi--Yau threefold has been checked by
B.~Szendr\H{o}i (see section 6.2 in \cite{szen}).
In this special case, the cohomology action induced by the
Fourier--Mukai transformation
has been written down by P. Aspinwall (section 6.1 in \cite{PA})
based on string theoretic arguments (apparently the formula is a reinterpretation of
previous physics results obtained in \cite{SDS}). Miraculously (at least
from the present author's point of view),
it can be checked that Aspinwall's formula is in complete agreement with the results
of this work.
\end{remark}

\subsection{Proof of the theorem}

It is not surprising that in order to prove theorem \ref{thm} we need
to compute the kernels $\tER \star \tEL,$
$\tEL \star \tER.$ As an intermediate step we will study the kernels
$\tRE \star \tLE$ and $\tLE \star \tRE.$ To that end, 
a few more facts about kernels (or correspondences)
in $\Der (Y),$ $Y = E \times_Z E$ are needed.

An object ${\Ga}_Y \in \Der (Y)$ of finite Tor--dimension 
defines an exact functor from $\Der (E)$
to $\Der (E)$ by the formula
\begin{equation*}
\Phi_{\tG_Y}(\;?):=\Rb q_{2_*}
(\tG_Y \Ltensor q_1^*(\;?)).
\end{equation*}
Note that we work under the important assumption that the morphism
$q: E \to Z$ is flat. We can define the composition of two kernels
$\tG_Y$ and $\tF_Y$ by the analog of the formula (\ref{def:comp})
\begin{equation*}
\tG_Y \star \tF_Y :=\Rb q_{13_*} \big( 
q_{12}^* (\tF_Y) \big)
\Ltensor
q_{23}^* (\tG_Y),
\end{equation*}
where the projection maps are flat morphisms from $E \times_Z E \times_Z E$
to $E \times_Z E.$
All the properties mentioned above for the composition of kernels in
$\Der (X \times X)$ continue to hold in this case, most remarkably the
associativity (due to the flatness of $q: E \to Z$).

Note also that we made the choice to denote the operation of composition of
kernels in $Y$ by the same symbol as before. It will be clear from the specific
context which composition is meant in a particular formula.

The functors $i_*$ and $\Lb i^*$ determined by the embedding
$i: E \hookrightarrow X$ can be expressed as exact functors determined
by the kernel $\Gamma_* (\Oa_E) \in \Der (E \times X)$
viewed as a correspondence from $E$ to $X,$ or from $X$ to $E,$ respectively
($\Gamma : E \hookrightarrow E \times X$ is the graph embedding). It 
is a nice and
easy exercise to check, using (\ref{eq:assoc}), that any kernel
$l_*(\Fa) \in \Der  (X \times X)$ can be decomposed as
\begin{equation}\label{eq:decomp}
l_*(\Fa) \cong \Gamma_* (\Oa_E) \star \Fa \star \Gamma_* (\Oa_E),
\end{equation}
($l : E \times E \hookrightarrow X \times X$ is the canonical embedding).

\begin{lemma} \label{lemma:decomp}

If $\tF_Y, \tG_Y \in \Der (Y),$ then
\begin{equation*}
k_*(\tG_Y) \star k_*(\tF_Y) \cong k_*(\tG_Y \star \tF_Y),
\end{equation*}
where $k : Y \to E \times E$ is the canonical embedding.

\end{lemma}

\begin{proof} (of the lemma) 

We need to show that
\begin{equation*}
k_* \Rb q_{13_*} \big( q_{12}^* (\tF_Y)
\Ltensor
q_{23}^* (\tG_Y)\big) \cong
\Rb r_{13_*} \big( r_{12}^* (k_*(\tF_Y)) \Ltensor
r_{23}^* (k_*(\tG_Y)) 
\big).
\end{equation*}
Consider the fiber square (with $r_{12}$ flat)
\begin{equation}
\begin{split}
\xymatrix@C=4pc{
E_1 \times_Z E_2 \times E_3 \ar[r]^-{m_{12}} \ar[d]_{(k, Id_{E_3})}
& E_1 \times_Z E_2 \ar[d]^k\\
E_1 \times E_2 \times E_3 \ar[r]^-{r_{12}}                    & E_1 \times E_2
}
\end{split}
\end{equation}
where we have used subscripts to distinguish between different copies of $E.$
We have that $r_{12}^*  k_* \cong (\Lb k^*, Id^*_{E_3})  m_{12}^*.$
By the projection formula, we can write
\begin{equation}\label{eq:ffff}
\begin{split}
&\Rb r_{13_*} \big( r_{12}^* (k_*(\tF_Y)) \Ltensor
r_{23}^* (k_*(\tG_Y)) 
\big) \cong \\
&\cong
\Rb r_{13_*} (k_*, Id_{E_3*}) \big( m_{12}^* (\tF_Y) \Ltensor
(\Lb k^*, Id^*_{E_3})
r_{23}^* (k_*(\tG_Y)) \big)
\end{split}
\end{equation}
By using another fiber square (with $r_{23} \circ (k, Id_{E_3})$ flat)
\begin{equation}
\begin{split}
\xymatrix@C=4pc{
E_1 \times_Z E_2 \times_Z E_3 \ar[r]^-{q_{23}} \ar[d]_{k_{12}}
& E_2 \times_Z E_3 \ar[d]^k\\
E_1 \times_Z E_2 \times E_3 \ar[r]^-{r_{23} \circ (k, Id_{E_3})}
& E_2 \times E_3
}
\end{split}
\end{equation}
Hence $ (\Lb k^*, Id^*_{E_3}) r_{23}^* k_* \cong k_{12_*} q^*_{23},$ and again
by the projection formula we obtain that the right hand side of (\ref{eq:ffff})
is isomorphic to
\begin{equation*}
\begin{split}
&\Rb r_{13_*} \big( m_{12}^* (\tF_Y) \Ltensor
k_*, Id_{E_3*}) \big( k_{12_*} q^*_{23} (\tG_Y) 
\big) \cong \\
&\cong \Rb r_{13_*} ( k_{12}^* m_{12}^* (\tF_Y) \Ltensor
k_*, Id_{E_3*}) k_{12_*} \big( q^*_{23} (\tG_Y) 
\big)  \cong k_* \Rb q_{13_*} \big( q_{12}^* (\tF_Y) \Ltensor
q_{23}^* (\tG_Y) \big),
\end{split}
\end{equation*}
since $r_{13} \circ (k, Id_{E_3}) \circ k_{12} = k \circ q_{13}.$
This ends the proof of the lemma.
\end{proof}

\begin{proposition}\label{prop:star}
If $\Ea \in \Der (E)$ is $EZ$--spherical, there exist 
distinguished triangles 
\begin{eqnarray}
&\tLE \longrightarrow \tRE \star \tLE \longrightarrow
\tRE \longrightarrow  \tLE \; [1], &\\
&\tLE \longrightarrow \tLE \star \tRE
\longrightarrow
\tRE \longrightarrow  \tLE \; [1].& 
\end{eqnarray}
\end{proposition}

\begin{proof} Write 
$$
\tRE= j_* \tRE_Y, \tLE= j_* \tLE_Y,
$$
with 
$$
\tRE_Y :=  q_1^* (\Du_E \Ea \Ltensor i^! \Oa_X)
\; \Ltensor \; q_2^* (\Ea), \ 
\tLE_Y := q_1^* (\Du_E \Ea \Ltensor q^! \Oa_Z)
\; \Ltensor \; q_2^* (\Ea).
$$
By (\ref{eq:decomp}), lemma \ref{lemma:decomp}
and the associativity of the composition of correspondences we can write that
\begin{equation*}
\begin{split}
\tRE \star \tLE &= j_*(\tRE_Y) \star j_*(\tLE_Y) \cong \\
&\cong \Gamma_* (\Oa_E) \star k_*(\tRE_Y)  \star
\big(\Gamma_* (\Oa_E) \star \Gamma_* (\Oa_E) \big) \star k_*(\tLE_Y) \star
\Gamma_* (\Oa_E) \\
\end{split}
\end{equation*}
We can now interpret the last line of the previous formula as
an exact functor $T : \Der (E \times E) \to \Der (X \times X)$
given by
$$
T(\Ka):= \Gamma_* (\Oa_E) \star k_*(\tRE_Y)  \star \Ka
\star k_*(\tLE_Y) \star \Gamma_* (\Oa_E).
$$
Lemma \ref{lemma:sga} shows that there exists a Postnikov 
system of the type (\ref{diagr:postnikov})
with 
$a=-d, b=0, \Ka = \Gamma_* (\Oa_E) \star \Gamma_* (\Oa_E),$
and $\Ha^{-c} (\Ka) \cong \Delta_* \Lambda^{c} \nu^\vee,$ 
for $0 \leq c \leq d,$
where $\Delta : E \hookrightarrow E \times E$ is the 
diagonal embedding.

In order to apply part $i)$ of lemma \ref{lemma:usus}, we need
to examine the action of the functor $T$ on the cohomology
sheaves $\Ha^{-c}.$ We have that
\begin{equation*}
\begin{split}
T(\Ha^{-c}[c])&= \Gamma_* (\Oa_E) \star k_*(\tRE_Y)  \star
\Delta_* \Lambda^{c} \nu^\vee [c]
\star k_*(\tLE_Y) \star
\Gamma_* (\Oa_E) \cong \\
&\cong \Gamma_* \Oa_E \star k_* \big( \tRE_Y \star 
\Delta_{Y*} \Lambda^{c} \nu^\vee [c]
\star \tLE_Y  \big) \star
\Gamma_* \Oa_E,
\end{split}
\end{equation*}
where we have used lemma \ref{lemma:decomp}, and the fact that
$\Delta = k \circ \Delta_Y$ (see diagram (\ref{diagr:zz})).

We have that
\begin{equation*}
(\tRE_Y \star \Delta_{Y*} \Lambda^{c} \nu^\vee [c] ) \star \tLE_Y 
\cong \Rb q_{13_*} \Big( q_{12}^* \tLE_Y \Ltensor
q_{23}^* \tRE_Y \star
\Delta_{Y*} \Lambda^{c} \nu^\vee [c] \Big).
\end{equation*}
But 
\begin{equation*}
\begin{split}
q_{23}^* \tRE_Y \star
\Delta_{Y*} \Lambda^{c} \nu^\vee [c] &\cong 
(q^{123}_2)^* (\Du_E \Ea \Ltensor (\Lambda^{c} \nu^\vee [c] \Ltensor
i^! \Oa_X) ) \Ltensor (q^{123}_3)^* (\Ea) \cong \\
q_{12}^* \tLE_Y &\cong  (q^{123}_1)^* (\Du_E \Ea \Ltensor q^! \Oa_Z) \Ltensor
(q^{123}_2)^* (\Ea).
\end{split}
\end{equation*}
Since 
\begin{equation}\label{eq:ghghp}
\Lambda^{c} \nu^\vee [c] \Ltensor i^! \Oa_X \cong \Lambda^{d-c} \nu [-d+c],
\end{equation}
after regrouping the terms, we obtain that
\begin{equation*}
\begin{split}
&\tRE_Y \star \Delta_{Y*} \Lambda^{c} \nu^\vee [c] \star \tLE_Y \cong 
\Rb q_{13_*} \big( q_{13}^* \tLE_Y \Ltensor (q^{123}_3)^*
(\Rhh_E (\Ea, \Ea \Ltensor \Lambda^{d-c} \nu [-d+c] ) \big) \cong \\
&\cong 
\tLE_Y \Ltensor \Rb q_{13_*} (q^{123}_3)^* (\Rhh_E (\Ea, 
\Ea \Ltensor \Lambda^{d} \nu [-d+c] ) \cong \\
&\cong \tLE_Y \Ltensor t^*
\big( \Rb q_* \Rhh_E ( \Ea, \Ea \Ltensor \Lambda^{d-c} \nu [-d+c] ) \big).
\end{split}
\end{equation*}
Since $\Ea$ is $EZ$--spherical, the last line in the above formula is
zero unless $c=0,d.$ For $c=d,$ the last line is $\tLE_Y,$ and for
$c=0,$ it is, by remark \ref{rem:mmm},
\begin{equation}\label{eq:ana}
\tLE_Y \Ltensor t^* \theta [-d-k] \cong q_1^* (\Du_E \Ea \Ltensor q^! \Oa_Z
\Ltensor q^* \theta [-d-k]) \; \Ltensor \; q_2^* (\Ea) \cong \tRE_Y,
\end{equation}
since $q^! \Oa_Z \Ltensor q^* \theta [-d-k] \cong i^! \Oa_X,$ 
by the starting assumption (\ref{def:theta}) made on $\theta.$ 
Part $i)$ of lemma \ref{lemma:usus} implies then
that there exist a distinguished triangle
$$
T(\Ha^{-d}[d]) \cong \tLE \longrightarrow
T(\Gamma_* \Oa_E \star \Gamma_* \Oa_E)
= \tRE \star \tLE \longrightarrow
T(\Ha^{0}) \cong \tRE \longrightarrow \tLE [1].
$$
Note that the proof works also in the case $d=0.$ In that case, lemma 
\ref{lemma:usus} is not needed, and the above calculation
shows that 
$$
\tRE \star \tLE \cong T(\Oa_X) \cong \tLE_Y \Ltensor t^*
\big( \Rb q_* \Rhh_E ( \Ea, \Ea \big),
$$
and the distinguished triangle (\ref{eq:ttrr}) defining the $EZ$--condition
in this case finishes the argument.

The proof of the existence of the second distinguished triangle 
is very similar, and requires the reversals of the roles 
of $\tRE$ and $\tLE,$ 
the replacement of the observation (\ref{eq:ghghp}) by 
$$
\Lambda^{c} \nu^\vee [c] \Ltensor q^! \Oa_Z \cong \Lambda^{d-c} \nu [-d+c]
\Ltensor q^* \theta^{-1} [d+k],
$$
and the use of the required analog of (\ref{eq:ana}), namely
$$
\tRE_Y \Ltensor t^* \theta^{-1} [d+k] \cong q_1^* (\Du_E \Ea \Ltensor i^! \Oa_X
\Ltensor q^* \theta^{-1} [d+k]) \; \Ltensor \; q_2^* (\Ea) \cong \tLE_Y.
$$
\end{proof}

\begin{proposition} \label{prop:hlp}
If $\Ea \in \Der (E)$ is $EZ$--spherical, then
\begin{equation*}
\tEL \star \tRE \; [1] \cong \tLE \cong \tRE \star \tEL \; [1].
\end{equation*}
\end{proposition}

\begin{proof}
For the first isomorphism, note that 
the operation $\star$ with one argument fixed
is an exact functor of triangulated categories, so
definition \ref{def:lr} provides the triangle
\begin{equation}\label{eq:disre2}
\tEL \star \tRE \longrightarrow
\tRE 
\overset{l_{\Ea}^{\#}}\longrightarrow
\tLE \star \tRE \longrightarrow
\tEL \star \tRE \; [1].
\end{equation}
We also consider the distinguished triangle of proposition (\ref{prop:star})
\begin{equation}\label{eq:disre1}
\tLE \longrightarrow \tLE \star \tRE
{\overset{g_{\Ea}^{\#} } \longrightarrow}
\tRE \longrightarrow  \tLE \; [1].
\end{equation}

We start with a lemma. 

\begin{lemma} 
The morphism $g_{\Ea}^{\#} \circ l_{\Ea}^{\#} \in \Hom_{X \times X} (\tRE, \tRE)$ 
is an isomorphism.
\end{lemma}

\begin{proof} (of the lemma) Since $\Ea$ is 
$EZ$--spherical,  proposition \ref{prop:uuuse} $iii)$ implies that 
$\Hom_{X \times X} (\tRE, \tRE) \cong k,$ so it is enough to show that
the morphism $g_{\Ea}^{\#} \circ l_{\Ea}^{\#}$ is non--zero. 
In fact, we will show that their composition 
\begin{equation*}
\tRE  \overset{l_{\Ea}^{\#}}\longrightarrow
\tLE \star \tRE \overset{g_{\Ea}^{\#}}\longrightarrow
\tRE,
\end{equation*}
induces group homomorphisms of Hom groups 
(all isomorphic to $\kb$) by proposition \ref{prop:uuuse} $i),iii)$)
\begin{equation*}
\Hom_{X \times X} (\tRE, \Delta_* \Oa_X ) 
\overset{g_{\Ea}^{\%}}\longrightarrow
\Hom_{X \times X} (\tLE \star \tRE, \Delta_* \Oa_X ) 
\overset{l_{\Ea}^{\%}}\longrightarrow
\Hom_{X \times X} (\tRE, \Delta_* \Oa_X ),
\end{equation*}
that are isomorphisms, where, proposition \ref{prop:uuuse} $i)$
shows that
$$
\Hom_{X \times X} (\tRE, \Delta_* \Oa_X ) \cong \kb,
$$
and,
propositions \ref{prop:adjoint}
and \ref{prop:uuuse} $iii),$ imply that
\begin{equation*}
\Hom_{X \times X} (\tLE \star \tRE, \Delta_* \Oa_X ) \cong 
\Hom_{X \times X} (\tLE, \tLE) \cong \kb.
\end{equation*}

In order to show that the composition $l_{\Ea}^{\#} \circ g_{\Ea}^{\#}$
is non-zero, we will ``probe'' it with the help of the 
cohomological
functor $\Hom_{X \times X} (\; ?, \Delta_* \Oa_X ).$ First, we apply this 
functor to
the distinguished triangle (\ref{eq:disre1}), and we 
look at a piece of the resulting long exact sequence of groups. 
\begin{equation*}
\begin{split}
\ldots &\to
\Hom_{X \times X} (\tLE [1], \Delta_* \Oa_X ) \to \\
&\to \Hom_{X \times X} (\tRE, \Delta_* \Oa_X ) 
\overset{g_{\Ea}^{\%}}\longrightarrow
\Hom_{X \times X} (\tLE \star \tRE, \Delta_* \Oa_X ) 
\to \ldots.
\end{split}
\end{equation*}
But by proposition \ref{prop:uuuse} $ii),$ 
the leftmost group is zero, hence the 
morphism $g_{\Ea}^{\%}$ is in fact an isomorphism 
\begin{equation*}
k \cong \Hom_{X \times X} (\tRE, \Delta_* \Oa_X ) \sito
\Hom_{X \times X} (\tLE \star \tRE, \Delta_* \Oa_X ) \cong \kb.
\end{equation*}

We now apply the same cohomological functor 
to the distinguished triangle (\ref{eq:disre2}) and investigate
a piece of the corresponding long exact sequence of groups.
\begin{equation*}
\begin{split}
\ldots &\to
\Hom_{X \times X} (\tEL \star \tRE [1], \Delta_* \Oa_X ) \to \\
&\to 
\Hom_{X \times X} (\tLE \star \tRE, \Delta_* \Oa_X ) 
\overset{l_{\Ea}^{\%}}\longrightarrow
\Hom_{X \times X} (\tRE, \Delta_* \Oa_X ) \to \ldots.
\end{split}
\end{equation*}
Proposition \ref{prop:adjoint} implies that
\begin{equation}\label{eq:www}
\Hom_{X \times X} (\tEL \star \tRE [1], \Delta_* \Oa_X ) \cong
\Hom_{X \times X} (\tEL, \tLE [-1])
\end{equation}
and proposition \ref{prop:corse} shows the 
the group is zero. Hence the 
morphism $l_{\Ea}^{\%}$ is an isomorphism 
\begin{equation*}
k \cong \Hom_{X \times X} (\tLE \star \tRE, \Delta_* \Oa_X ) \sito
\Hom_{X \times X} (\tRE, \Delta_* \Oa_X ) \cong \kb,
\end{equation*}
which finishes the proof of the lemma.
\end{proof}

The lemma shows that we can 
consider the following ``9--diagram'' (page 24 in \cite{BBD}) which is essentially
a version of the octahedron axiom in triangulated categories.

\begin{equation}
\begin{split}
\xymatrix{
{\tRE [1]} \ar[r]^{\cong} & \tRE [1]  \ar[r]  & 0  \ar[r]  & \tRE [2] \\
\tEL \star \tRE [1] \ar[r] \ar[u] & 0 \ar[r] \ar[u] & \wt{\mathcal U} \ar[u] \ar[r] &
\tEL \star \tRE [2] \ar[u]  \\
\tLE \star \tRE \ar[u] \ar[r]^{g_{\Ea}^{\#}}
  & \tRE \ar[r] \ar[u] & \tLE [1] \ar[r] \ar[u]
& \tLE \star \tRE [1] \ar[u] \\
\tRE \ar[r]^{\cong} \ar[u]^{l_{\Ea}^{\#}}
 & \tRE \ar[r] \ar[u]^{Id} & 0 \ar[r] \ar[u] & \tRE [1] \ar[u]
}
\end{split}
\end{equation}
The starting point is the commutative square located in the lower left corner. 
The ``9--diagram'' proposition
(prop. 1.1.11. in \cite{BBD}) shows that we can fill in the diagram, so
$\tEL \star \tRE [2] \cong
\wt{\mathcal U} \cong \tLE [1].$ This ends the proof of the first isomorphism
of the proposition.

The proof of the second one is very similar. Of course, we have to replace 
the distinguished triangles (\ref{eq:disre2}) and (\ref{eq:disre1})
by 
\begin{equation*}
\begin{split}
&\tRE  \star \tEL \longrightarrow \tRE 
{\overset{l_{\Ea}^{\prime\#}} \longrightarrow}
\tRE \star \tLE \longrightarrow
\tRE \star \tEL \; [1],\\
&  \tLE 
\longrightarrow \tRE \star \tLE  {\overset{g_{\Ea}^{\prime\#}} 
\longrightarrow}
\tRE \longrightarrow  \tLE \; [1],
\end{split}
\end{equation*}
and then show that the morphism $g_{\Ea}^{\prime\#} \circ l_{\Ea}^{\prime\#} \in 
\Hom_{X \times X} (\tRE, \tRE)$ is an isomorphism.

Everything works as above-- the only modification is that
the isomorphism (\ref{eq:www}) has to be replaced by 
\begin{equation*}
\Hom_{X \times X} (\tRE \star \tEL [1], \Delta_* \Oa_X ) \cong
\Hom_{X \times X} (\tEL, \tLE [-1]),
\end{equation*}
which is true due to proposition \ref{prop:adjoint}. 

\

We are now in the position to finish the proof of the theorem. The argument is
similar to the one used to prove the previous proposition. 
Consider the distinguished triangle
(\ref{eq:basiclr}) defining $\tEL$
\begin{equation}
\label{eq:basiclra}
\tEL \to \Delta_* \Oa_X
\overset{l_{\Ea}}\longrightarrow
\tLE \overset{u_{\Ea}}\longrightarrow
\tEL [1],
\end{equation}
with a choice of a (non--canonical) morphism $u_\Ea.$
Since 
$\Ea \in \Der(E)$ is $EZ$--spherical, 
we apply the 
cohomological functor $\Hom_{X \times X} (\tLE, \; ?)$ 
to this distinguished triangle and write a relevant 
piece of the associated long exact sequence: 
\begin{equation*}\begin{split}
\ldots &\to
\Hom_{X \times X} (\tLE, \Delta_* \Oa_X )  
\overset{l_{\Ea}^{\%}}\longrightarrow
\Hom_{X \times X} (\tLE, \tLE) \overset{u_{\Ea}^{\%}}\longrightarrow
\Hom_{X \times X} (\tLE, \tEL [1]) \to \\
& \to \Hom_{X \times X} (\tLE, \Delta_* \Oa_X [1] ) \to \ldots. 
\end{split}
\end{equation*}
By proposition \ref{prop:uuuse} $ii), iii),$ 
we know that
\begin{equation*}
\Hom_{X \times X} (\tLE, \Delta_* \Oa_X ) \cong 0, \; 
\Hom_{X \times X} (\tLE, \tLE) \cong \kb, \\
\end{equation*}
and by propositions \ref{prop:adjoint} $i),$ \ref{prop:hlp} and 
\ref{prop:uuuse} $i)$ (in this order),
we have that
\begin{equation}\label{eq:klkkk}
\Hom_{X \times X} (\tLE, \tEL [1]) \cong 
\Hom_{X \times X} ( \Delta_* \Oa_X, \tRE \star \tEL [1]) \cong 
\Hom_{X \times X} ( \Delta_* \Oa_X, \tLE ) \cong \kb.
\end{equation}

We conclude that 
the 
morphism $u_\Ea \in \Hom_{X \times X} (\tLE, \tEL \; [1])$ is a
generator of this group.

Let's now look at the distinguished triangle,
\begin{equation}\label{eq:tri100}
\tEL \to \tEL \star \tER \to \tEL \star \tRE [1] \cong \tLE 
\overset{v_{\Ea}}\longrightarrow
\tEL [1],
\end{equation}
obtained by applying the exact functor $\tEL \star (?)$ to 
the distinguished triangle (\ref{eq:basiclr}) defining $\tER.$
Due to proposition \ref{prop:hlp}, we have that 
$\tEL \star \tRE [1] \cong \tLE,$ so 
\begin{equation*}
\Hom_{X \times X} (\tEL \star \tRE [1], \tEL \; [1]) \cong 
\Hom_{X \times X} (\tLE, \tEL \; [1])
\cong \kb,
\end{equation*}
by (\ref{eq:klkkk}).

We claim that the morphism $v_{\Ea}$ is a 
generator of this group (i.e. non--zero).
Indeed, after applying the cohomological
functor
$\Hom_{X \times X} (\tLE, \; ?)$ to the distinguished 
triangle (\ref{eq:tri100}), we can write a relevant piece of the
resulting long exact sequence as follows:
\begin{equation*}
\ldots \to  
\Hom_{X \times X} (\tLE, \tEL \star \tER) \to
\Hom_{X \times X} (\tLE, \tEL \star \tRE [1] ) \overset{v_{\Ea}^{\%}}\longrightarrow
\Hom_{X \times X} (\tLE, \tEL [1]) \to 
\ldots.
\end{equation*}
Proposition \ref{prop:hlp} and 
proposition \ref{prop:uuuse} $iii)$ show that 
\begin{equation*}
\Hom_{X \times X} (\tLE, \tEL \star \tRE [1] ) \cong
\Hom_{X \times X} (\tLE, \tLE) ) \cong \kb.
\end{equation*}
Moreover, we can write that
\begin{equation*}
\begin{split}
&\Hom_{X \times X} (\tLE, \tEL \star \tER) \cong
\Hom_{X \times X} ( \Delta_* \Oa_X, (\tRE \star \tEL[1]) \star \tER [-1]) \cong \\
&\cong \Hom_{X \times X} ( \Delta_* \Oa_X, \tLE \star \tER [-1]) \cong 
\Hom_{X \times X} (\tRE, \tER [-1]), 
\end{split}
\end{equation*}
where the first and the last isomorphisms follow from proposition 
\ref{prop:adjoint}, and the middle one from proposition 
\ref{prop:hlp}. Proposition \ref{prop:corse} 
says that the last group is zero,
which proves indeed that $v_{\Ea}$ is non--zero.

In other words, we have shown that the morphisms $u_{\Ea}$ and 
$v_{\Ea}$ coincide up to an isomorphism. We can then look
at the diagram 

\begin{equation}
\begin{split}
\xymatrix{
\tEL \ar[r] \ar[d]^{\cong} & \Delta_* \Oa_X 
\ar[r] \ar@{.>}[d] & \tLE \ar[r]^{u_{\Ea}} \ar[d]^{\cong} & \tEL [1] \ar[d]^{\cong} \\
\tEL \ar[r] & \tEL \star \tER \ar[r] & \tEL \star \tRE [1] 
\ar[r]^{v_{\Ea}} & \tEL [1]
}
\end{split}
\end{equation}
We have just argued that the right hand square is commutative. The
axiom TR3 of a triangulated category implies the existence
of the dotted morphism. Another well known property of
triangulated categories (see, for example, corollary 4, page 242 
in \cite{gelmanin}) shows that it is also an isomorphism, which 
concludes the proof of theorem \ref{thm}.
\end{proof}

\section{Applications}\label{sec:appl}

In this section, we sample some geometric geometric situations
 in which $EZ$--spherical objects arise. The list has no claims of being exhaustive.

\begin{example} \label{exam:st}
The case $Z={\rm Spec} (\kb)$ brings nothing new, namely
an object $\Ea \in \Der(E)$ is $EZ$--spherical if and only if
$i_* \Ea \in \Der(X)$ is spherical, and
the functors $\Phi_{\tER}$  and $\Phi_{\tEL}$
are those studied in \cite{ST}.
This follows immediately by adjunction for the pair $(i_*, i^!)$ and
by recalling that, in this case, the functors $\Rb^{\bullet}q_* \Rhh_E$
are nothing else but the ${\rm Ext}^\bullet$ groups on $E.$
\end{example}

\begin{example} Let $E$ be a smooth proper divisor in a smooth
quasi--projective
variety $X,$ and $q=Id : E \to Z=E.$ 
Then $\theta \cong \omega_E
\Ltensor \Lb i^* \omega_X^{-1} \cong \omega_{E/X} \cong \nu.$ 
In this case $Y= E \times_Z E  \cong E,$ and any invertible 
sheaf $\Ea$ is $EZ$--spherical, since
\begin{equation*}
\Rb q_* \Rhh_E (\Ea, \Ea) \cong \Oa_Z
\end{equation*}
The object $\tEL$ in $\Der (X \times X)$ is 
$\Delta_* \Oa(-E),$ and the associated 
Fourier--Mukai functor is the usual 
\begin{equation*}
\Phi_{\tEL} (?) = \Oa(-E) \Ltensor (?).
\end{equation*}
\end{example}

\begin{example} Assume that $E=X.$ This is the case of the so--called
fiberwise Fourier--Mukai functors \cite{BrMa}, \cite{thomas},
\cite{ACRY}. The typical example is the case of a flat projective 
Calabi--Yau ($\omega_X \cong \Oa_X$)
fibration $q: X \to Z$ with a generic Calabi--Yau fiber $F$ of dimension $k.$
Consider an object $\Ea \in \Der(X)$ flat over the base $Z.$

Assume that the restriction $\Ea_z$ is a spherical object (in the sense of
\cite{ST}) for any generic fiber $F_z,$ that the dimensions of
${\rm Ext}^i (\Ea_z, \Ea_z)\; (z \in Z)$ are constant functions on $Z$
for all $i,$
and that $q_* \Rhh_X (\Ea, \Ea) \cong \Oa_Z.$
The Grauert--Grothendieck 
theorem (\cite{hart2} III.12.9)
implies that
\begin{equation*}
\Rb^i q_* \Rhh_X (\Ea, \Ea) =0, \ \text{for} \ 0 < i < k.
\end{equation*}
In this case $\theta = \omega_Z,$ and 
since $q^! \Oa_Z \cong q^* \omega_Z [k],$ the duality
theorem (\ref{eq:dual}) (or relative Serre duality)
gives that $\Ea$ is $EZ$--spherical.

When the generic fiber $F$ has the property that $h^{i,0} (F)=0,$
for $0 < i < k,$ any invertible sheaf on $X$ will be an $EZ$--spherical
object. For a flat elliptic fibration with a section $\sigma : Z \to X,$
the sheaf $\Oa_{\sigma(Z)}$ is also $EZ$--spherical.

\end{example}

\

The next example shows one possible way of relating our results
to the theory of exceptional objects \cite{rudakov}.

\begin{example} Assume that $E$ is a proper smooth divisor in $X$
such that $E \cong F \times Z$ ($F$ comes from either fiber, or Fano),
$q : E \to Z,$ $\dim F=k.$ Assume that $\theta$ is an invertible sheaf on $Z$ 
such that $q^! \theta \cong \omega_{E/X} [k]$
Let $\Ea^\prime$ be an exceptional object in $\Der(F),$ i.e.
${\rm Hom}_F (\Ea', \Ea')=k,$ and ${\rm Ext}^c_F (\Ea', \Ea')=0,$ for
$c > 0.$ Then the pull-back of $\Ea'$ to $E$ by the natural projection
is an $EZ$--spherical object in $\Der(E).$
\end{example}

\begin{example} \label{exam:mori}
Let's also discuss in more detail the example that provided the inspiration
for this work. Of course, it can generalized to cover more general
situations in birational geometry; however, 
this example captures the essential
features. We study the $EZ$--spherical objects that arise when one
considers a smooth Calabi--Yau complete intersection $X$ in a toric variety $\wt X$
and a toric  elementary contraction $\wt q : \wt X \to \wt Y$
in the sense of Mori theory (as studied by M. Reid in \cite{reid}). As
mentioned in the introduction, such transformations provide geometric
models for the physical phase transitions.
Let $\wt E$ represent the loci where $\wt q$ is not an
isomorphism (the exceptional
locus) and $\wt Z:= \wt q (\wt E).$ Corollary 2.6 in \cite{reid} shows that
$\wt E$ is a complete intersection of toric Cartier
divisors $D_1, \ldots, D_d,$ and the restriction of $\wt q$ to $\wt E$ is a
flat morphism whose fibers are weighted projective spaces.

Let $E$ denote the exceptional locus of the restriction of the
contraction to $X$ and $q: E \to Z$ the corresponding morphism.
We make the assumption that
\begin{equation*}
E = X \; \cap \; \wt E,
\end{equation*}
and that $q$ is a flat morphism with the fibers $F$ (regular) projective spaces
of dimension $k.$ For more
on the geometrical details of the situation, the reader may want to consult
section 4.4 in \cite{rph}.

According to proposition 2.7 in \cite{reid}
we have that
\begin{equation}\label{eq:curve}
D_c \cdot C_\sigma < 0,
\end{equation}
for all $c, 1 \leq c \leq d,$ and for any curve $C_\sigma$ in the class of
the contraction, therefore
any rational curve contained in a fiber $F \cong \PP^k.$

On the other hand, the
adjunction formula implies that
\begin{equation}\label{eq:adj00}
\omega_F = \Oa (D_1 + \ldots + D_d)_{|_F}.
\end{equation}
The only way that (\ref{eq:curve}) and (\ref{eq:adj00}) could hold at same time
is if all the line bundles of the type
$\Oa (D_{c_1} + \ldots + D_{c_m})_{|_F}$ on $F \cong \PP^k$ with $0 < m < d$
are negative of the form $\Oa (-l)$ with $0 < l < k+1.$

We see that in this case
\begin{equation*}
\nu_{|_E} \cong \big( \Oa (D_1) \oplus \ldots \oplus \Oa(D_d) \big)_{|_E},
\end{equation*}
and
\begin{equation*}
{\rm H}^l (F,  \Lambda^c \nu_{|_F})=0, \ \text{for} \ 0 < l < k+1, 0 < c < d.
\end{equation*}
Definition \ref{def:EZ} and
the Grauert--Grothendieck 
theorem (\cite{hart2} III.12.9)
imply that any invertible sheaf $\La$ on $E$ is indeed an $EZ$--spherical
object, since $\Rhh_E (\La, \La)= \Oa_E$ and $\Rb q_* \Oa_E = \Oa_Z.$

Note that there is nothing special about the Calabi--Yau condition. To put
the $EZ$--machinery to work, we only need to know that the restriction to $E$
of the canonical bundle of $X$ is also the pull--back of an invertible
sheaf on $Z.$
It should also be noted that, in the course of analyzing this example, the
existence of the contraction of $X$ is not needed. Presumably, there could
be examples with $E$ a complete intersection of Cartier divisors
$D_1, \ldots, D_d$ in $X$ and $q : E \to Z$ a flat fibration
with a Fano fiber $F.$ In that case, we would only need to impose as
an assumption the analog of (\ref{eq:curve}) in order to have that
any invertible sheaf on $E$ is $EZ$--spherical.
\end{example}

\begin{remark} \label{rem:normal}
An $EZ$--spherical object $\Ea$ has good ``portability''
properties. Indeed, definition \ref{def:EZ} has a local character with
respect to the embedding $E \hookrightarrow X.$ Therefore,
if $\Ea$ is $EZ$--spherical with respect to a
configuration described by a diagram such as (\ref{diag:gen}), it will
continue to be
$EZ$--spherical if the variety $X$ is replaced (analytically)
by the ``local'' variety given by the total space of the normal
bundle $N_{E/X},$ while $q : E \to Z$ is left unchanged.
In fact, $\Ea$ will remain $EZ$--spherical if $X$ is replaced by any
other smooth variety $X'$ such that $N_{E/X'} \cong N_{E/X}.$
The explicit passage from $X$ to $X'$ can be realized
by using the formal completion
of $X$ (or $X'$) along $E$ (\cite{hart2} II.9), or by the so--called
``deformation to the normal cone'' (\cite{fulton} chap. 5).
\end{remark}



\begin{thebibliography}{10}

\bibitem{ACRY}
B.~Andreas, G.~Curio, D.~Hernandez Ruiperez and S.-T.~Yau, {\em Fourier--Mukai
transform and mirror symmetry for D--branes on elliptic Calabi--Yau},
preprint (2000), {\tt math.AG/0012196}.

\bibitem{PA}
P.~S. Aspinwall, {\em Some navigation rules for D--brane monodromy},
J. High Energy Phys. {\bf 2} 2001, {\tt hep-th/0102198}.

\bibitem{AspDon}
P.~S. Aspinwall and R.~Y. Donagi, {\em The heterotic string, the tangent
bundle, and derived categories}, Adv. Theor. Math. Phys. {\bf 2} (1998),
1041--1074, {\tt hep-th/9806094}.

\bibitem{AGM1}
P.~S. Aspinwall, B.~R. Greene and D.~R. Morrison, {\em {C}alabi--{Y}au
moduli space, mirror manifolds and spacetime topology change in
string theory},
  Nuclear Phys. B {\bf 416} (1994), 414--480, {\tt hep-th/9309097}.

\bibitem{ahk}
P.~S. Aspinwall, R.~P. Horja, R.~L. Karp, {\em 
Massless D-Branes on Calabi-Yau Threefolds and Monodromy},
preprint (2002), {\tt hep-th/0209161}.

\bibitem{BBD}
A.~A. Beilinson, J.~N. Bernstein and P.~Deligne, {\em Faisceaux pervers},
Ast\'erisque {\bf 100} (1982).

\bibitem{sga6}
P.~Berthelot, A.~Grothendieck and L.~Illusie, 
{\em Th\'eorie des intersections et th\'eor\`eme de Riemann-Roch,} 
S\'eminaire de G\'eom\'etrie Alg\'egrique du Bois-Marie 1966--1967 (SGA 6).
Springer-Verlag, Berlin-New York (1971).

\bibitem{BO}
A.~Bondal and D.~Orlov, {\em Semiorthogonal decompositions
for algebraic varieties}, preprint (1995), {\tt alg-geom/9506012}.

\bibitem{Bridge}
T.~Bridgeland, {\em Equivalences of triangulated categories and
Fourier--Mukai transforms}, Bull. London Math. Society {\bf 31} (1999),
25--34, {\tt math.AG/9809114}.

\bibitem{Bridge1}
T.~Bridgeland, {\em Flops and derived categories}, 
Invent. Math. {\bf 147} (2002), 613--632, {\tt math.AG/0009053}.

\bibitem{Bridgenew}
T.~Bridgeland, {\em Stability conditions on triangulated categories},
preprint (2002), {\tt math.AG/0212237}. 

\bibitem{BKR}
T.~Bridgeland, A.~King and M.~Reid, {\em Mukai implies McKay: the McKay
correspondence as an equivalence of derived categories},  
J. Amer. Math. Soc. {\bf 14} (2001), 535--554, 
{\tt math.AG/9908027}.

\bibitem{BrMa}
T.~Bridgeland and A.~Maciocia, {\em Fourier-Mukai transforms for K3
and elliptic fibrations},  J. Amer. Math. Soc. {\bf 14} (2001), 535--554, 
{\tt math.AG/9808022}.

\bibitem{CKYZ}
T.-M. Chiang, A.~Klemm, S.-T. Yau and E.~Zaslow, {\em Local
mirror symmetry: calculations and interpretations},
Adv. Theor. Math. Phys. {\bf 3} (1999), 495--565, {\tt hep-th/9903053}.

\bibitem{BConrad}
B.~Conrad, {\em Grothendieck Duality and Base Change},
Lect. Notes Math. 1750,
Springer--Verlag, Berlin Heidelberg New York (2000).


\bibitem{CK}
D.~A. Cox and S.~Katz, {\em Mirror symmetry and algebraic geometry},
Math. Surveys and Monographs, Vol. 68, American Math. Soc., 1999.

\bibitem{alishii}
A.~Craw, A.~Ishii, {\em Flops of G-Hilb and equivalences of derived 
categories by variation of GIT quotient}, 
preprint (2002), {\tt  math.AG/0211360}.

\bibitem{deligne1}
P.~Deligne, {\em Cohomologie a support propre et construction du
foncteur $f^!$}, Appendix to \cite{hart1}.

\bibitem{Douglas1}
M.~R. Douglas, {\em D-branes, categories and $N=1$ supersymmetry},
J. Math. Phys. {\bf 42} (2001), 2818--2843,
{\tt  hep-th/0011017}.

\bibitem{fulton}
W.~Fulton, {\em Intersection theory}, Second edition,
Springer--Verlag, Berlin (1998).

\bibitem{GKZ2}
I.~M. Gelfand, M.~M. Kapranov and A.~V. Zelevinsky,
{\em Discriminants, resultants and multidimensional determinants},
Birkh\"auser Boston, 1994.

\bibitem{gelmanin}
S.~I. Gelfand and Y.~I. Manin,
{\em Methods of homological algebra}, Translated from the 1988 Russian
original, Springer--Verlag, Berlin (1996).

\bibitem{glo}
V.~Golyshev, V.~Lunts and D.~Orlov, {\em Mirror symmetry for abelian
varieties}, J. Algebraic Geom. {\bf 10} (2001), 433--496,
{\tt math.AG/9812003}.

\bibitem{hart1}
R.~Hartshorne, {\em Residues and duality}, Lect. Notes Math. 20,
Springer--Verlag, Berlin Heidelberg New York (1966).

\bibitem{hart2}
R.~Hartshorne, {\em Algebraic geometry}, Springer--Verlag,
Berlin Heidelberg New York (1977).

\bibitem{rph}
R.~P. Horja, {\em Hypergeometric functions and mirror symmetry in
toric varieties}, preprint (1999), {\tt math.AG/9912109}.

\bibitem{KatzMayrVafa}
S.~Katz, P.~Mayr and C.~Vafa, {\em Mirror symmetry and exact solution
of 4D N=2 gauge theories I}, Adv. Theor. Math. Phys. {\bf 1} (1998),
53--114, {\tt hep-th/9706110}.

\bibitem{SDS}
S.~Katz, D.~R. Morrison and M.~R. Plesser, {\em Enhanced gauge
symmetry in type II string theory}, Nuclear Phys. B {\bf 477} (1996),
105--140, {\tt hep-th/9601108}.

\bibitem{kleiman}
S.~L. Kleiman, {\em Relative duality for quasicoherent sheaves},
Compositio Math. {\bf 41} (1980), 39--60.

\bibitem{Kont1}
M.~Kontsevich, {\em Homological algebra of mirror symmetry}, Proc.
Internat. Congr. Math. Z\"urich 1994 (S.~D. Chatterji, ed.), vol.~1,
Birkh\"auser Verlag, Basel, Boston, Berlin, 1995, 120--139,
{\tt alg-geom/9411018}.

\bibitem{Kont2}
M.~Kontsevich, Lecture at Rutgers University, November 11, 1996,
(unpublished).

\bibitem{KontSoib}
M.~Kontsevich and Y.~Soibelman, {\em Homological mirror symmetry and
torus fibrations}, 
Symplectic geometry and mirror symmetry (Seoul, 2000), 203--263, 
World Sci. Publishing, River Edge, NJ, {\tt math.SG/0011041}.

\bibitem{LoLu}
E.~Looijenga and V.~L. Lunts, {\em A Lie algebra attached to a projective
variety}, Inv. Math. {\bf 129} (1997), 361--412, {\tt
alg-geom/9604014}.

\bibitem{maninecm}
Yu.~I. Manin, {\em Moduli, Motives, Mirrors},
Plenary talk at the 3rd ECM, Barcelona, July 10-14, 2000,
{\tt math.AG/0005144}.

\bibitem{Mukai1}
S.~Mukai, {\em Duality between $D(X)$ and $D(\tilde{X})$
with its application to Picard sheaves}, Nagoya Math. J. {\bf 81}
(1981), 153--175.

\bibitem{O}
D.~Orlov, {\em Equivalences of derived categories and $K3$
surfaces}, Algebraic geometry 7, J. Math. Sci. (New York) {\bf 84} (1997),
1361--1381, {\tt alg-geom/9606006}.

\bibitem{reid}
M.~Reid, {\em Decomposition of toric morphisms}, Arithmetic and geometry,
vol.~II, Birkh\"auser Boston (1983), 395--418.

\bibitem{rudakov}
A.~N. Rudakov Et Al, {\em Helices and Vector Bundles,} Seminaire
Rudakov, London Math. Soc. Lecture Note Series {\bf 148},
Cambridge Univ. Press (1990).

\bibitem{seidel}
P.~Seidel, {\em Lagrangian two-spheres can be symplectically knotted},
J. Diff. Geom. {\bf 52} (1999), 145--171, {\tt math.DG/9803083}.

\bibitem{ST}
P.~Seidel and R.~Thomas, {\em Braid group actions on derived
categories of coherent sheaves}, 
Duke Math. J. {\bf 108} (2001), 37--108, {\tt math.AG/0001043}.

\bibitem{szen}
B.~Szendr\H{o}i, {\em Diffeomorphisms and families of Fourier-Mukai transforms in
mirror symmetry},
Applications of algebraic geometry to coding theory, physics and computation 
(Eilat, 2001), 317--337, Kluwer Acad. Publ., Dordrecht, 2001,  
{\tt math.AG/0103137}.

\bibitem{balazs1}
B.~Szendr\H{o}i, \hspace{0.1cm}
{\em Artin group actions on derived categories of threefolds}, \hspace{0.1cm}
preprint,  (2002), {\tt math.AG/0210121}.

\bibitem{balazs2}
B.~Szendr\H{o}i, \hspace{0.2cm}
{\em Enhanced gauge symmetry and braid group actions}, \hspace{0.2cm}
preprint \hspace{0.15cm} (2002), {\tt  math.AG/0210122}.

\bibitem{thomas}
R.~P. Thomas, {\em Mirror symmetry and actions of braid groups
on derived categories}, Proceedings of the Harvard Winter School
on Mirror Symmetry, International Press (1999), {\tt math.AG/0001044}.

\bibitem{verdier0}
J.-L.~Verdier, {\em Base change for twisted inverse image of coherent sheaves},
Algebraic Geometry (Internat. Colloq., Tata Inst. Fund. Res., Bombay, 1968)
393--408, Oxford Univ. Press, London.

\bibitem{verdier}
J.-L.~Verdier, {\em Des cat\'egories d\'eriv\'ees des cat\'egories
ab\'eliennes,} With a preface by Luc Illusie (G.~Maltsiniotis, ed.),
Ast\'erisque, {\bf 239} (1996).

\bibitem{Witten1}
E.~Witten, {\em Phases of $N=2$ theories in two dimensions},
Nucl. Physics B {\bf 403} (1993), 159--222, {\tt hep-th/9301042}.

\end{thebibliography}
\end{document}